\documentclass{amsart}

\usepackage{amscd,amssymb}
\def\Z{{\mathbb Z}}
\def\Q{{\mathbb Q}}
\def\R{{\mathbb R}}
\def\C{{\mathbb C}}
\def\L{{\mathbb L}}
\def\N{{\mathbb N}}

\def\V{{\mathbb V}}
\def\W{{\mathbb W}}

\def\Ga{{\mathbb G}_a}					
\def\Schur{{\mathbb S}}					

\def\A{{\mathcal A}}
\def\G{{\mathcal G}}
\def\H{{\mathcal H}}
\def\J{{\mathcal J}}
\def\M{{\mathcal M}}
\def\O{{\mathcal O}}
\def\U{{\mathcal U}}
\def\T{{\mathcal T}}

\def\cL{{\mathcal L}}
\def\cC{{\mathcal C}}
\def\cP{{\mathcal P}}
\def\cR{{\mathcal R}}

\def\d{{\mathfrak d}}
\def\f{{\mathfrak f}}
\def\g{{\mathfrak g}}
\def\h{{\mathfrak h}}
\def\n{{\mathfrak n}}
\def\o{{\mathfrak o}}
\def\p{{\mathfrak p}}
\def\r{{\mathfrak r}}
\def\s{{\mathfrak s}}
\def\t{{\mathfrak t}}
\def\u{{\mathfrak u}}

\def\sp{\s\p}						

\def\rhotilde{\tilde{\rho}}
\def\omegatilde{\widetilde{\omega}}
\def\Jbar{\overline{\J}}
\def\Text{\widehat{T}}

\def\dot{{\bullet}}
\def\Efin{E_{{\mathrm fin}}}		

\def\blank{\phantom{x}}
\def\sup#1{^{(#1)}}
\def\A#1{A_k^{(#1)}}
\def\B#1{B_k^{(#1)}}
\def\bil#1#2{\langle #1,#2 \rangle}
\newcommand\comptensor{\operatorname{\widehat{\otimes}}}

\def\LiE{{\textsf{LiE}}}
\def\rr{\raggedright}

\newcommand\im{\operatorname{im}}		
\newcommand\id{\operatorname{id}}
\newcommand\ad{\operatorname{ad}}

\newcommand\Hom{\operatorname{Hom}}
\newcommand\End{\operatorname{End}}
\newcommand\Ext{\operatorname{Ext}}
\newcommand\Aut{\operatorname{Aut}}
\newcommand\Der{\operatorname{Der}}

\newcommand\Jac{\operatorname{Jac}}

\newcommand\Gr{\operatorname{Gr}}

\def\Hcts{H_{\mathrm{cts}}}

\newtheorem{theorem}{Theorem}[section]
\newtheorem{lemma}[theorem]{Lemma}
\newtheorem{proposition}[theorem]{Proposition}
\newtheorem{corollary}[theorem]{Corollary}

\theoremstyle{definition}
\newtheorem{definition}[theorem]{Definition}

\theoremstyle{remark}
\newtheorem{remark}[theorem]{Remark}

\newtheorem{question}[theorem]{Question}

\begin{document}

\title{Infinitesimal Presentations of the Torelli Groups}

\author{Richard Hain}

\address{Department of Mathematics\\ Duke University\\
Durham, NC 27708-0320}

\email{hain@math.duke.edu}

\date{July 18, 1996}

\thanks{This work was supported in part by grants from the National
Science Foundation.}


\maketitle

\tableofcontents

\section{Introduction}

By a theorem of Malcev \cite{malcev}, every torsion free nilpotent
group can be imbedded canonically as a discrete, cocompact subgroup
of a real nilpotent Lie group. One can therefore associate to a
finitely generated group $\pi$ a tower of nilpotent Lie groups
\begin{equation}\label{tower}
\cdots \to G_3 \to G_2 \to G_1 = H_1(\pi,\R)
\end{equation}
by taking $G_k$ to be the nilpotent Lie group associated to the
maximal torsion free quotient of $\pi$ of length $k$. Since each
nilpotent Lie group is simply connected, the tower (\ref{tower})
is determined by the corresponding tower
$$
\cdots \to \g_3 \to \g_2 \to \g_1 = H_1(\pi,\R)
$$
of nilpotent Lie algebras. The inverse limit $\g$ of this
tower is a pronilpotent Lie algebra, called the
{\it Malcev Lie algebra associated to $\pi$.} This Lie algebra
has the property that the graded Lie algebra $\Gr \g$
associated with its lower central series is isomorphic to
$\left(\Gr \pi\right) \otimes \R$, where $\Gr \pi$ is the
graded $\Z$-Lie algebra associated to the filtration of $\pi$ by
its lower central series.

In this paper we give a presentation of the Malcev Lie algebra
$\t_{g,r}^n$ associated to the Torelli group $T_{g,r}^n$ for all
$g\ge 6$. Since each Torelli group injects into its unipotent
completion (at least when $n+r>0$), the corresponding Malcev
Lie algebra should contain significant information about the group.

Recall that the mapping class group $\Gamma_{g,r}^n$ \label{group_def}
is defined as follows. Fix a compact orientable surface $S$ of
genus $g$, together with $n+r$ distinct points
\begin{equation}\label{points}
x_1,\dots,x_n; y_1,\dots,y_r
\end{equation}
of $S$ and $r$ non-zero tangent vectors $v_1,\dots,v_r$, where
$v_j$ is tangent to $S$ at $y_j$. The group $\Gamma_{g,r}^n$
is the group of isotopy classes of orientation preserving
diffeomorphisms of $S$ that fix each of the points (\ref{points})
and each of the tangent vectors $v_j$.%
\footnote{One can replace each tangent vector in the  definition
by a boundary component --- with this change, the diffeomorphisms
are required to be the identity on each boundary component.}
The Torelli group $T_{g,r}^n$ \label{torelli_def}is defined to be
the kernel of the natural homomorphism
\begin{equation}\label{rho}
\Gamma_{g,r}^n \to \Aut H_1(S,\Z).
\end{equation}
Observe that the classical pure braid group $P_n$ is $T_{0,1}^n$.

Our presentation of $\t_{g,r}^n$ generalizes the well-known
presentation of $\p_n$, the Malcev Lie algebra of the pure braid
group $P_n$, which is
of importance in the theory of Vassiliev invariants (cf.\
\cite{kontsevich:vass,cartier:knots,bar-natan}) and was first
written down by Kohno \cite{kohno:braids}. Denote the free Lie
algebra generated by indeterminates $X_1,\dots, X_m$ by
$\L(X_1,\dots,X_m)$, and that generated by a vector space $V$
by \label{lie_def} $\L(V)$. Then $\p_n$ is the completion of the
graded Lie algebra
$$
\L(X_{ij} : ij\text{ is a two element subset of }\{1,\dots,n\})/R
$$
where $R$ is the ideal generated by the quadratic relations
\begin{align*}\label{braid_relns}
[X_{ij},X_{kl}]&\text{ when $i,j,k$ and $l$ are distinct;}\cr
[X_{ij},X_{ik} + X_{jk}]& \text{ when $i,j$ and $k$ are distinct.}
\end{align*}
The property that $\p_n$ is the completion of the associated graded
Lie algebra $\Gr \p_n$ does not hold for the generic group, but does
hold for all Torelli groups as we shall see.

It is easiest to first state the result for the absolute Torelli
group, $T_g := T_{g,0}^0$. It follows from Dennis Johnson's  computation
of the first homology of $T_g$ \cite{johnson:h1} that each graded
quotient of the lower central series of $\t_g$ is a representation
of the algebraic group $Sp_g$. We will give a presentation of
$\Gr \t_g$ in the category of representations of $Sp_g$. Chose a set
$\lambda_1,\dots, \lambda_g$ of fundamental weights of $Sp_g$.
Denote the representation of $Sp_g$ with highest weight
$\lambda = \sum n_j \lambda_j$ by $V(\lambda)$. Johnson's
fundamental computation is that there is a natural $Sp_g(\Z)$
equivariant isomorphism between $H_1(T_g,\Q)$ and $V(\lambda_3)$.

For all $g \ge 3$, the representation $\Lambda^2 V(\lambda_3)$
contains a unique copy of $V(2\lambda_2) + V(0)$. Denote the
$Sp_g$ invariant complement of this by $R_g$. Since the quadratic
part of the free Lie algebra $\L(V)$ is $\Lambda^2 V$, we can
view $R_g$ as being a subspace of the quadratic elements of
$\L(V(\lambda_3))$.

\begin{theorem}
For all $g\neq 2$, $\t_g$ is the completion of its
associated graded $\Gr\t_g$. When $g\ge 6$, this has presentation
$$
\Gr \t_g = \L(V(\lambda_3))/(R_g),
$$
where $R_g$ is the set of quadratic relations defined above.
When $3 \le g < 6$, the relations in $\Gr \t_g$ are generated by
the quadratic relations $R_g$, and possibly some cubic relations.
\end{theorem}

In fact, we will show that in genus 3 there are no quadratic relations
and the cubic relations contain a copy of $V(\lambda_3)$.

Dennis Johnson has proved that $T_g$ is finitely generated for
all $g\ge 3$, but it is not known for any $g\ge 3$ whether or not
$T_g$ is finitely presented. Geoff Mess \cite{mess} proved that
$T_2$ is a countably generated free group. (Note that when
$g = 0,1$, $T_g$ is trivial.)

\begin{corollary}
For all $g\neq 2$, and for all $r,n \ge 0$, $t_{g,r}^n$ is
finitely presented. \qed
\end{corollary}

In the decorated case, we have the extension
\begin{equation}\label{extension}
1 \to \pi_{g,r}^n \to T_{g,r}^n \to T_g \to 1,
\end{equation}
where $\pi_{g,r}^n$ denotes the fundamental group of the
configuration space of $n$ points and $r$ tangent vectors
in $S$. After applying the Malcev Lie algebra functor, we
obtain an extension
$$
0 \to \p_{g,r}^n \to \t_{g,r}^n \to \t_g \to 0,
$$
where $\p_{g,r}^n$ denotes the Malcev Lie algebra of
$\pi_{g,r}^n$.

\begin{theorem}
For all $g \ge 0$, and all $r,n \ge 0$, the Lie algebra
$\p_{g,r}^n$ is the completion of its associated graded
$\Gr \p_{g,r}^n$. The associated graded has a presentation
with only quadratic relations.
\end{theorem}

The explicit presentation is given in Section \ref{braids2}.
In order to give the presentation for $\t_{g,r}^n$, we
prove that (\ref{extension}) remains exact after taking
graded quotients. Thus, in order to give a presentation
of $\t_{g,r}^n$, it suffices to determine the map
$$
[\blank,\blank] :
\left(\Gr^1 \t_g \otimes \Gr^1 \p_{g,r}^n\right)
\oplus \left(\Gr^1 \t_g \otimes \Gr^1 \t_g\right)
\to \Gr^2 \t_{g,r}^n.
$$
determined by the bracket.
We do this in Section \ref{decorated} to obtain the presentation
of $\t_{g,r}^n$ in general.

Our results complement, and sometimes overlap with, the beautiful
work \cite{morita:casson,morita:cocycles,morita:trace,morita:conj} of
Shigeyuki
Morita who began the study of the ``higher Johnson homomorphisms''
studied in this paper. Our main theorem allows us to answer several
questions about Torelli groups, and to prove a conjecture of Morita.
These and other applications are discussed in Section~\ref{applications}.

Another feature of the classical case is the existence of a
canonical universal integrable connection. Denote the
classifying space
$$
\C^n - \left\{(z_1,\dots,z_n) :
\text{ the $z_i$ are not distinct}\right\}
$$
of $P_n$ by $X_n$. Denote the complex of global meromorphic
$k$-forms on a
complex manifold $Y$ by $\Omega^k(Y)$. The universal integrable
connection on $X_n$ is given by the $\p_g$ valued 1-form
$$
\sum_{ij} d\log(z_i - z_j)\, X_{ij} \in \Omega^1(X_n) \otimes \p_n.
$$
It plays a central role in the theory of Vassiliev invariants
(cf. \cite{kohno:KZ}, \cite{cartier:knots}, \cite{kassel}.)
We are able to prove that there is a canonical universal
connection form with ``scalar curvature'' for each $T_{g,r}^n$,
provided $g\neq 2$, although, to date, we have not been able to
give an explicit formula for it. The universal connection is
discussed in Section \ref{applications}.

The basic approach in this paper is to use Hodge
theory. The main technical theorem of the paper is:

\begin{theorem}
Suppose that $g\neq 2$ and that $r,n \ge 0$. For each choice of
a complex structure on the decorated reference surface
$$
\left(S;x_1,\dots,x_n;y_1,\dots,y_r;v_1,\dots,v_r\right)
$$
there is a mixed Hodge structure on $\t_{g,r}^n$ for
which the bracket is a morphism of mixed Hodge structures.
\end{theorem}

This mixed Hodge structure is canonical once one fixes
an isomorphism of $\Gamma_{g,r}^n$ with the (orbifold) fundamental
group of the moduli space of smooth projective curves of genus
$g$ with $n$ marked points, and $r$ non-zero tangent vectors.
The theorem is proved using the mixed Hodge structure on
the completion of the mapping class group $\Gamma_{g,r}^n$
relative to the homomorphism $\Gamma_{g,r}^n \to Sp_g(\Q)$
induced by (\ref{rho}), the existence of
which follows from \cite{hain:derham}. The mixed Hodge
structure on the Lie algebra $\u_{g,r}^n$ of the prounipotent
radical of the relative completion is lifted to $\t_{g,r}^n$
using two results from \cite{hain:comp}. The first states that
we have a central extension
$$
0 \to \Ga \to \t_{g,r}^n \to \u_{g,r}^n \to 0
$$
when $g\ge 3$.
The second gives an explicit relationship between this extension
and the algebraic 1-cycle $C - C^-$ in the jacobian of an
algebraic curve $C$. The theory of relative completion is
reviewed in Section \ref{rel_comp}.

When $g \ge 3$, the weight filtration of $\t_{g,r}^n$ is its
lower central series. The fact that the weight graded functors
are exact on the category of mixed Hodge structures then
allows the reduction to the associated graded with impunity
when studying $\t_{g,r}^n$, $\u_{g,r}^n$, $\p_{g,r}^n$ and
maps between them.

In order to bound the degrees of relations in $\t_g$ by $N$,
we need to know that if $\V$ is a variation of Hodge structure
of weight $n$ over $\M_g$ (the moduli space of curves) that comes
from a rational representation of $Sp_g$, then the weights on
$$
H^2(\M_g,\V)
$$
are bounded between $2+n$ and $n+N$ --- see Section \ref{cts_coho_tor}.
There is no {\it a priori} uniform bound on the weights of $H^k(X,\V)$,
where $X$ is a smooth variety and $\V$ is a variation of Hodge structure
over $X$ of weight $l$, as there is in the case of $\Q$ coefficients
where the weights are bounded between $k$ and $2k$. For example, if
$\Gamma$ is a finite index subgroup of $SL_2(\Z)$ and $X$ the quotient
of the upper half plane by $\Gamma$, then the non-trivial weights on
$H^1(X,S^n\V)$ are $n+1$ and $2n+2$ for infinitely many $n$, as can
be seen from results in \cite{zucker}. Here $\V$ denotes the
fundamental representation of $SL_2$ viewed as a variation of
Hodge structure over $X$ of weight 1 and $S^n\V$ its $n$th symmetric
power. Thus, one of the main technical ingredients in the paper is the
result of Kabanov \cite{kabanov} (see also \cite{kabanov:purity} which
states that one can take $N$ to be 2 when $g\ge 6$, and 3 when
$3 \le g < 6$.

The existence of the mixed Hodge structure on the Malcev Lie algebra
associated to the Torelli group was obtained several years ago. The
quadratic relations (proved in Section \ref{quadratic_relns}) were
derived in  \cite{hain:letter}. Subsequently Morita (unpublished)
proved that when the genus is sufficiently large there are no cubic
or quartic relations in $\t_g$. Kabanov's purity theorem allows us to
avoid  Morita's involved computations and to show there are no
higher order relations.
\medskip

\noindent{\it Acknowledgements.} It is a pleasure to thank all those
with whom I have had useful discussions, especially A.~Borel, P.~Deligne,
Alexander Kabanov, Eduard Looijenga, Shigeyuki Morita, and Steven Zucker.
I would also like to thank Hiroaki Nakamura for his numerous comments
on the manuscript.
I would like to thank the Institute for Advanced Study, the Institut des
Hautes \'Etudes Scientifiques, the Institut Henri Poincar\'e, and the
Universit\"at Essen, each of which supported me during my sabbatical
during which this paper was written.

\section{Braid Groups in Positive Genus}
\label{braids1}

Throughout this section $g$ will be positive. Suppose that
$S$ is a compact oriented surface of genus $g$, and that
$r$ and $n$ are integers $\ge 0$. The configuration space of
$m\ge 1$ points on $S$ is
$$
F^m(S) = S^m - \Delta,
$$
where $\Delta$ is the union of the various diagonals $x_i = x_j$.
Denote the tangent bundle of $S$ by $TS$, and the bundle of non-%
zero tangent vectors by $V$. The pullback of $V$ to $F^m(S)$ along
the $j$th projection $p_j : F^m(S) \to S$ will be denoted by $V_j$.
For a subset $A$ of $\{1,\dots,m\}$ denote the fibered product
of the $V_j$, where $j \in A$, by $V_A$.

The {\it configuration space \label{config_def} $F_{g,r}^n$ of $n$
points and $r$ non-zero tangent vectors of $S$} is defined to be the
total space of the bundle
$$
V_A \to F^{r+n}(S)
$$
where $A = \{n+1,\dots,n+r\}$. Fix a base point $f_o$ of $F_{g,r}^n$.
Define \label{fund_def}
$$
\pi_{g,r}^n = \pi_1(F_{g,r}^n,f_o).
$$
When $r=0$ this is just the group of pure braids with $n$ strings
on the surface $S$. In general, this group can be thought of as the
group of pure braids on $S$ with $r+n$ strings where $r$ of the
strings are framed. It is a standard fact that the space $F_{g,r}^n$
is an Eilenberg-MacLane space of type $K(\pi,1)$ \cite[\S 1.2]{birman}.

In contrast with the genus 0 case, we have:

\begin{proposition}\label{h1_braid}
For each $g\ge 0$, there is a short exact sequence
$$
0 \to \left(\Z/(2g-2)\Z\right)^r \to H_1(\pi_{g,r}^n,\Z)
\stackrel{p}{\to} H_1(S^{n+r},\Z) \to 0,
$$
where $p$ is induced by the natural map $F_{g,r}^n \to S^{n+r}$.
\end{proposition}

\begin{proof}
We first consider the case when $r=0$. In this case $F_{g,r}^n$
is $S^n - \Delta$. The divisor $\Delta$ is the
union of the diagonals $\Delta_{ij}$ where the $i$th and $j$th
point of $S^n$ are equal. We therefore have a Gysin
sequence
$$
\cdots \to H_2(S^n) \stackrel{\gamma}{\to} \bigoplus_{i<j}\Z
\stackrel{t}{\to} H_1(S^n - \Delta) \to H_1(S^n) \to 0.
$$
The map $\gamma$ takes a cycle $z$ to the element of
$\oplus_{i<j}\Z$ whose $ij$th term is the intersection
number $z\cdot\Delta_{ij}$. The map $t$ takes the
generator of the $ij$th factor to the homology class of a
small circle which winds about about $\Delta_{ij}$ in the
positive direction.

Fix a base point $x_o$ of $S$.
For $u\in H_k(S)$ and $i\in \{1,\dots,n\}$, denote by
$u^i$ the element
$$
x_o \times \dots \times x_o \times \stackrel{i}{u}
\times x_o \times \dots \times x_o
$$
of $H_k(S^n,\Z)$, where $u$ is placed in the $i$th factor.
For elements $u$ and $v$ of $H_1(S)$ and $i,j\in\{1,\dots,n\}$,
denote the element
$$
x_o \times \dots \times x_o \times \stackrel{i}{u}
\times x_o \times \dots  \times x_o \times \stackrel{j}{v}
\times x_o \times \dots \times x_o
$$
of $H_2(S,\Z)$ by $u^i \times v^j$, where $u$ is placed in
the $i$th factor and $v$ in the $j$th.

By choosing representatives of $u$ and $v$ which do not
pass through the base point, one sees immediately that
$$
\gamma : u^i \times u^j \mapsto (u\cdot v)\Delta_{ij}
$$
from which it follows that $\gamma$ is surjective and that
$t$ is trivial. This proves the result when $r=0$.

Observe that
$$
\gamma : S^i \mapsto \sum_{j\neq i} \Delta_{ij}.
$$
If $a$ and $b$ are elements of $H_1(S,\Z)$ with intersection
number 1, then
$$
S^i - \sum_{j\neq i} a^i \times b^j
$$
is in the kernel of $\gamma$, and therefore lifts to
an element $\sigma_i$ of $H_2(S^n-\Delta,\Z)$.

We prove the general case by induction on $r$. Our
inductive hypothesis is that the result has been proven
for $F_{g,s}^m$ when $s<r$, and that there are classes
$$
\sigma_1,\dots, \sigma_m \in H_2(F_{g,s}^m,\Z)
$$
whose images under the the maps
$$
{p_j}_\ast : H_2(F_{g,s}^m,\Z) \to H_2(S,\Z)
$$
induced by the various projections $p_j : F_{g,s}^m \to S$
satisfy $p_j(\sigma_i) = \delta_{ij}[S]$. We have proved this
when $r=0$.

Suppose that $r>0$. We have the projection
$$
q: F_{g,r}^n \to F_{g,r-1}^{n+1},
$$
which replaces the first tangent vector by its anchor point.
This is a principal $\C^\ast$ bundle. It fits into a commutative
square
$$
\begin{CD}
F_{g,r}^n 		@>>> V 	\cr
@VqVV 		  	@VVV 	\cr
F_{g,r-1}^{n+1} @>>> S 	\cr
\end{CD}
$$
One therefore has a map
$$
\begin{CD}
H_2(F_{g,r-1}^{n+1}) @>\gamma_F>> H_0(F_{g,r}^n) @>>>
H_1(F_{g,r}^n) @>>> H_1(F_{g,r-1}^{n+1}) @>>> 0 \cr
			@VVV						@VVV
@VVV
		@VVV \cr
H_2(S)			@>\gamma_S>> H_0(S)		 @>>>
H_1(V)		 @>>> H_1(S)				@>>> 0 \cr
\end{CD}
$$
of Gysin sequences. The map $\gamma_S$ is simply multiplication
by the euler characteristic.
Since $\sigma_{n+1} \mapsto [S]$, we see that $\gamma_F$ takes
$\sigma_{n+1}$ to $2 - 2g$. The result follows by induction.
\end{proof}

\begin{remark}\label{error}
This result (with $r=0$) can also be proved by considering the
natural fibrations $F_g^{n+1} \to F_g^n$ obtained by forgetting
the last point. The fiber is an $n$ punctured copy of $S$, and
its homology therefore fits into an exact sequence
$$
0 \to \Z^n/\text{diagonal} \to
H_1(S - n\text{ points}) \to H_1(S) \to 0.
$$
One has to be a little careful as the monodromy action is not
trivial. A simple geometric argument shows that the monodromy
acts trivially on the kernel and quotient in the sequence above,
and therefore is given by a homomorphism
$$
\pi_g^n \to \Hom(H_1(S),\Z^n/\text{diagonal})
\cong H^1(S^n)/\text{diagonal}.
$$
By induction on $n$, $H_1(\pi_g^n)$ is isomorphic to $H_1(S^n)$.
Since the monodromy is abelian, it factors through the quotient
map
$$
\pi_g^n \to H_1(S)^{\oplus n}.
$$
A straightforward geometric argument shows that the action of the
latter is given by the map
$$
H_1(S)^{\oplus n} \stackrel{PD^{\oplus n}}{\longrightarrow}
H^1(S)^{\oplus n}/\text{diagonal},
$$
where $PD$ denotes Poincar\'e duality.
The coinvariants are therefore given by
$$
H_0(F_g^n,H_1(\text{fiber})) = H_1(S).
$$
An elementary spectral sequence argument completes the inductive
step.

Kohno and Oda \cite[p.~208]{kohno-oda} use this method, but their
result contradicts ours as they mistakenly assume that the
monodromy representation is trivial.
\end{remark}

\section{Relative Completion of Mapping Class Groups}
\label{rel_comp}

In this section we recall the main theorem of \cite{hain:comp}
which  makes precise the relationship between the Malcev completion
of $T_{g,r}^n$ and the unipotent radical of the relative Malcev
completion of $\Gamma_{g,r}^n$. We first recall the definition
of relative Malcev completion, which is due to Deligne. A reference
for this material is \cite[\S\S 2--4]{hain:comp}.

Suppose that $\Gamma$ is a discrete group, $S$ a reductive linear
algebraic group over a field $F$ of characteristic zero, and that
$\rho : \Gamma \to S(F)$ is a representation whose image is Zariski
dense. The {\it Malcev completion of $\Gamma$ over $F$ relative to
$\rho$} is a homomorphism $\rhotilde : \Gamma \to \G$ of $\Gamma$ into
a proalgebraic group $\G$, defined over $F$, which is an extension
$$
1 \to \U \to \G \stackrel{p}{\to} S \to 1
$$
of $S$ by a prounipotent group $\U$ such that the diagram
$$
\begin{CD}
\Gamma @>{\rho}>> S \cr
@V{\rhotilde}VV @| \cr
\G @>>p> S
\end{CD}
$$
commutes. It is characterized by a universal mapping property: If
$G$ is a linear (pro)algebraic group over $F$ which is an extension
$$
1 \to U \to G \to S \to 1
$$
of $S$ by a (pro)unipotent group, and if $\tau : \Gamma \to G$
is a homomorphism whose composition with $G \to S$ is $\rho$, then
there is a unique homomorphism $\G \to G$ such that the diagram
$$
\begin{CD}
\Gamma @>{\rhotilde}>> \G \cr
@V{\tau}VV		@VpVV\cr
G @>>> S
\end{CD}
$$
commutes.

When $S$ is the trivial group, the relative completion of $\Gamma$
coincides with the classical Malcev (or unipotent) completion of
$\Gamma$.

Suppose that $K/F$ is an extension of fields of characteristic
zero. When $S$ is defined over $F$ and $\rho : \Gamma \to S(F)$, one can
ask if the $K$-form of the completion of $\Gamma$ relative to $\rho$ is
obtained from the $F$-form by extension of scalars. If this is the case
for all such field extensions, we will say that the relative completion
of $\Gamma$ relative to $\rho$ can be defined over $F$.

The action of the mapping class group on $S$ preserves the intersection
pairing $q : H_1(S,\Z)^{\otimes 2} \to Z$. We therefore have a
homomorphism
\begin{equation}\label{map}
\rho : \Gamma_{g,r}^n \to \Aut (H_1(S,\Z),q) \cong Sp_g(\Z).
\end{equation}
For a positive integer $l$, we define the {\it level $l$ subgroup}
$\Gamma_{g,r}^n[l]$ \label{level_def} to be the kernel of the induced
map
$$
\Gamma_{g,r}^n \to \Aut (H_1(S,\Z/l\Z),q) \cong Sp_g(\Z/l\Z).
$$
Here we interpret $Sp_g(\Z/l\Z)$ as the trivial group when $l=1$.

\begin{theorem}\label{rational}
For all $g\ge 3$ and all $l \ge 1$, the completion of the mapping
class group $\Gamma_{g,r}^n[l]$ relative to the homomorphism
$\rho : \Gamma_{g,r}^n[l] \to Sp_g(\Q)$ induced by (\ref{map})
is defined over $\Q$. \qed
\end{theorem}

This result was proved in \cite[(4.14)]{hain:comp} under the
assumption that $g \ge 8$ and that $l=1$. That the stronger
result is true follows from the strengthening \cite{borel:improved}
of Borel's stability theorem \cite{borel:triv,borel:twisted} for the
symplectic group, stated below, which ensures that
the hypothesis \cite[(4.10)]{hain:comp} is satisfied when $l\ge 1$
and $g \ge 3$.

\begin{theorem}\label{imp_borel}
Suppose that $V$ is an irreducible rational representation of
the algebraic group $Sp_g$ and that $\Gamma$ is a finite index subgroup
of $Sp_g(\Z)$. If $k < g$, then $H^k(\Gamma,V)$ vanishes when $V$ is
non-trivial, and agrees with the stable cohomology of $Sp_g(\Z)$ when
$V$ is the trivial representation. \qed
\end{theorem}

Denote the completion of $\Gamma_{g,r}^n$ relative to $\rho$ by
$\rhotilde : \Gamma_{g,r}^n \to \G_{g,r}^n$. \label{comp_def} Denote the
prounipotent radical of $\G_{g,r}^n$ by \label{ugp_def} $\U_{g,r}^n$, and
its Lie algebra by \label{ulie_def} $\u_{g,r}^n$.

\begin{proposition}\label{level}
If $g\ge 3$, then for all $l\ge 1$, the composite
$$
\Gamma_{g,r}^n[l] \hookrightarrow \Gamma_{g,r}^n \to \G_{g,r}^n
$$
is the completion of $\Gamma_{g,r}^n[l]$ relative to the restriction
of $\rho$ to $\Gamma_{g,r}^n[l]$.
\end{proposition}

\begin{proof}
This follows directly from results in \cite[\S4]{hain:comp} as
we shall explain. Denote the relative completion of $\Gamma_{g,r}^n[l]$
by $\G_{g,r}^n[l]$ and its prounipotent radical
by $\U_{g,r}^n[l]$. There is a natural map
$\U_{g,r}^n[l] \to \U_{g,r}^n$, the surjectivity of which follows from
(\ref{imp_borel}) and \cite[(4.6)]{hain:comp}. Injectivity follows
directly from (\ref{imp_borel}) and  \cite[(4.13)]{hain:comp}.
\end{proof}

We have an extension
$$
1 \to \U_{g,r}^n \to \G_{g,r}^n \to Sp_g \to 1
$$
of proalgebraic groups over $\Q$. The homomorphism $\rhotilde$
induces a map $T_{g,r}^n \to \U_{g,r}^n$. Denote the classical
Malcev completion of \label{comptor_def} $T_{g,r}^n$ by $\T_{g,r}^n$,
and its Lie algebra by \label{lietor_def} $\t_{g,r}^n$.
Since $\U_{g,r}^n$ is prounipotent, $\rhotilde$ induces a homomorphism
$$
\theta : \T_{g,r}^n \to \U_{g,r}^n
$$
of prounipotent groups.

The following theorem is the main result of \cite{hain:comp}.%
\footnote{There is a minor error in proof of the case ``$A_{g,r}^n$
implies $A_{h,r}^n$'' of the proof of \cite[(7.4)]{hain:comp}. It is
easily fixed.}
There it is proved for all $g\ge 8$, but in view of
(\ref{imp_borel}), it holds for all $g\ge 3$. (Cf.\ the third
footnote on page~76 of \cite{hain:comp}.)

\begin{theorem}\label{central_ext}
For all $g\ge 3$, the homomorphism $\theta$ is surjective and has
a one dimensional kernel isomorphic to $\Ga$ which is central in
$\T_{g,r}^n$ and is trivial as an $Sp_g(\Z)$ module. Moreover,
the extensions are all pulled back from that of $\T_g$; that is,
the diagram
$$
\begin{CD}
0 @>>> \Ga @>>> \T_{g,r}^n @>>> \U_{g,r}^n @>>> 1 \cr
@. @| @VVV @VVV \cr
0 @>>> \Ga @>>> \T_g @>>> \U_g @>>> 1
\end{CD}
$$
commutes. \qed
\end{theorem}

It is a standard fact that the sequence
$$
1 \to \pi_{g,r}^n \to \Gamma_{g,r}^n \to \Gamma_g \to 1
$$
is exact. Restricting to $T_g$, we obtain an extension
\begin{equation}
\label{exten}
1 \to \pi_{g,r}^n \to T_{g,r}^n \to T_g \to 1.
\end{equation}

\begin{proposition}\label{h1_tor}
If $g\ge 3$, then the extension (\ref{exten}) induces an
exact sequence
$$
0 \to H_1(\pi_{g,r}^n,\Q) \to H_1(T_{g,r}^n,\Q)
\to H_1(T_g,\Q) \to 0.
$$
\end{proposition}

\begin{proof}
It follows from (\ref{h1_braid}) that the natural map
$\pi_{g,r}^n \to \pi_g^{n+r}$ induces an isomorphism on
$H_1$ with rational coefficients. The corresponding
surjection $T_{g,r}^n \to T_g^{n+r}$ induces a map
$$
\begin{CD}
H_1(\pi_{g,r}^n,\Q) @>>> H_1(T_{g,r}^n,\Q) @>>> H_1(T_g,\Q) @>>> 0\cr
@VVV @VVV @| \cr
H_1(\pi_g^{n+r},\Q) @>>> H_1(T_g^{n+r},\Q) @>>> H_1(T_g,\Q) @>>> 0\cr
\end{CD}
$$
of exact sequences. Since the middle vertical map is a surjection,
and the two outside maps are isomorphisms, it follows that the middle
map is an isomorphism. It therefore suffices to prove the result when
$r=0$.

To prove this, we need to prove that the map
\begin{equation}\label{red}
H_1(\pi_g^n,\Q) \to H_1(T_g^n,\Q)
\end{equation}
is injective. We first remark that this is easily proved when
$n=1$ using the Johnson homomorphism
$$
\tau_g^1 : H_1(T_g^1) \to H_3(\Jac C).
$$
(Cf.\ \cite{johnson:def} and \cite[\S 3]{hain:normal}.)
The composition of $\tau_g^1$ with the map
$$
H_1(\pi_g^1,\Q) \to H_3(\Jac C,\Q)
$$
is easily seen to be
the map $H_1(C) \to H_3(\Jac C)$ which takes a class in
$H_1$ to its Pontrjagin product with the class of $C$ in
$H_2(\Jac C)$. Since this map is injective, it follows that
(\ref{red}) is injective when $n=1$.

The general case follows from this by considering the
maps $p_j : H_1(T_g^n) \to H_1(T_g^1)$ induced by the
$n$ forgetful maps $T_g^n \to T_g^1$.
\end{proof}

Define $\cP_{g,r}^n$ \label{comppi_def} to be the Malcev completion
of $\pi_{g,r}^n$ and \label{liepi_def} $\p_{g,r}^n$ to be the
corresponding Malcev Lie algebra.
Applying Malcev completion to (\ref{exten}) we obtain a sequence
$$
\p_{g,r}^n \to \t_{g,r}^n \to \t_g
$$
of Malcev Lie algebras.

\begin{proposition}\label{seq}
If $g\neq 2$, then the sequence
$$
0 \to \p_{g,r}^n \to \t_{g,r}^n \to \t_g \to 0
$$
associated to (\ref{exten}) is exact.
\end{proposition}

\begin{proof}
By \cite[(5.6)]{hain:cycles} it suffices to verify two conditions.
First, that $T_g$ acts  unipotently on $H^1(\pi_{g,r}^n,\Q)$; this
follows from (\ref{h1_braid}). The second condition there is satisfied
if, for example, the extension
$$
0 \to H_1(\pi_{g,r}^n,\Q) \to G \to T_g \to 1
$$
obtained by pushing (\ref{exten}) out along
$\pi_{g,r}^n \to H_1(\pi_{g,r}^n,\Q)$ is split. In our case this
follows from (\ref{h1_tor}) as the extension
above can be pulled back from the extension
$$
0 \to H_1(\pi_{g,r}^n,\Q) \to H_1(T_{g,r}^n,\Q)
\to H_1(T_g,\Q) \to 0
$$
which is split for trivial reasons.
\end{proof}

The standard homomorphism $\Gamma_{g,r}^n \to \Gamma_g$ induces a
homomorphism $\G_{g,r}^n \to \G_g$ of relative completions. The
inclusion $\cP_{g,r}^n \to \T_{g,r}^n$ induces a homomorphism
$\cP_{g,r}^n \to \U_{g,r}^n$. We therefore
have a sequence
$$
 \cP_{g,r}^n \to \G_{g,r}^n \to \G_g \to 1
$$
of proalgebraic groups.

\begin{lemma}\label{exactness}
If $g\ge 3$, then the sequence
$$
1 \to \cP_{g,r}^n \to \G_{g,r}^n \to \G_g \to 1
$$
is exact. In particular, the sequence
$$
0 \to \p_{g,r}^n \to \u_{g,r}^n \to \u_g \to 0
$$
is exact.
\end{lemma}

\begin{proof}
To prove the result, it suffices to prove that the sequence
$$
1 \to \cP_{g,r}^n \to \U_{g,r}^n \to \U_g \to 1
$$
is exact. But this follows immediately from
(\ref{central_ext}) and (\ref{seq}).
\end{proof}

Suppose that $\g$ is a finitely generated pronilpotent Lie algebra.
Denote the group of automorphisms of $\g$ by $\Aut\g$. Denote the
subgroup of $\Aut\g$ consisting of the elements that act trivially
on $H_1(\g)$ by $L^1\Aut \g$. Since the
action of an automorphism on the graded quotients of the lower central
series is determined by the action on the first graded quotient,
$\Aut \g$ is a proalgebraic group which is an extension
$$
1 \to L^1\Aut \g \to \Aut \g \to S \to 1
$$
of a closed subgroup $S$ of $\Aut H_1(\g)$ by the prounipotent group
consisting of those automorphisms of $\g$ that act trivially on the
graded quotients of the lower central series. Its Lie algebra is the
Lie algebra $\Der \g$ of derivations of $\g$. This is an extension
$$
0 \to L^1\Der \g \to \Der \g \to \s \to 0
$$
of the Lie algebra $\s$ of $S$ by the pronilpotent Lie algebra
of derivations of $\g$ that act trivially on the graded quotients
of the lower central series of $\g$.

If $\G$ is the prounipotent group corresponding to $\g$, then
$\Aut \G$ and $\Aut \g$ are isomorphic, as can be seen using the
Baker-Campbell-Hausdorff formula.

\begin{lemma}\label{rep}
For all $g\ge 0$ the natural action of $\Gamma_{g,r}^n$ on
$\pi_{g,r}^n$ induces a representation
$$
\G_{g,r}^n \to \Aut \p_{g,r}^n.
$$
\end{lemma}

\begin{proof}
Suppose that $g \ge 0$. The mapping class group $\Gamma_{g,r}^n$
acts on $\p_{g,r}^n$. We therefore have a homomorphism
\begin{equation}\label{act}
\Gamma_{g,r}^n \to \Aut \p_{g,r}^n.
\end{equation}
By (\ref{h1_braid}) we know that
$$
\Aut H_1(\p_{g,r}^n) = \Aut H_1(S)^{\oplus(n+r)}.
$$
There is a diagonal copy of $Sp_g$ contained in this group,
and it is easy to see that this is the Zariski closure of the
image of $\Gamma_{g,r}^n$ in $\Aut H_1(\p_{g,r}^n)$. It follows
that the Zariski closure of the image of (\ref{act}) is an
extension of this diagonal $Sp_g$ by a prounipotent group. Since
the homomorphism from $\Gamma_{g,r}^n$ to this copy of $Sp_g$ is
the standard representation, the universal mapping property of the
relative completion implies that (\ref{act}) induces a
homomorphism $\G_{g,r}^n \to \Aut\p_{g,r}^n$.
\end{proof}

\begin{remark}\label{sl2}
When $g=1$, the results (\ref{level}) and (\ref{central_ext})
are false.
That (\ref{level}) and (\ref{central_ext}) fail can be deduced
from \cite[(10.3)]{hain:derham}, a special case of which states that
there is a  natural isomorphism
$$
H^1(\M_1[l],S^n\V) \cong
\left(H^1(\U_{1}[l])\otimes S^n\V\right)^{SL_2}.
$$
Here $\M_1[l]$ denotes the moduli space of elliptic curves with a
level $l$ structure, $\V$ denotes the variation of Hodge structure
over $\M_1[l]$ of weight 1 corresponding to $H^1$ of the universal
elliptic curve, and $S^n \V$ denotes its $n$th symmetric power.
Since the level $l$ congruence subgroup of $SL_2(\Z)$ is free when
$g\ge 4$, it follows by an Euler characteristic argument that
$H^1(\M_1[l],S^k\V)$ is non-zero whenever $g\ge 4$. Since $\T_1$
is trivial, it cannot surject onto $\U_1[l]$. Since the rank $r_l$
of the level $l$ subgroup of $SL_2(\Z)$ depends on $l$, and since
$$
\dim H^1(\M_1[l],S^n\V) = (r_l-1)\dim S^k\V,
$$
it follows that the rank of the $S^n\V$ isotypical part of
$H_1(\U_1[l])$ depends on $l$. So (\ref{level}) does not hold.
\end{remark}

\section{Mixed Hodge Structures on Torelli Groups}

Denote by $\M_{g,r}^n[l]$ \label{mod_def} the moduli space of
ordered $(n+r+1)$-tuples
$$
(C;x_1,\dots,x_n;v_1,\dots,v_r)
$$
where $C$ is a smooth complex projective curve with a level
$l$ structure, the $x_j$ are distinct points of $C$, and the
$v_j$ are non-zero holomorphic tangent vectors of $C$ which
are anchored at $r$ distinct points of $C$ which are also distinct
from the $x_j$. We shall omit the $l$ when it is 1, and $r$ and
$n$ when they are zero. So, for example, $\M_g$ denotes the moduli
space of smooth projective curves of genus $g$.

For each point $x$ of $\M_{g,r}^n[l]$, there is a natural isomorphism
of $\Gamma_{g,r}^n[l]$ with the (orbifold) fundamental group
$\pi_1(\M_{g,r}^n[l],x)$ of $\M_{g,r}^n[l]$. We will denote the
latter by $\Gamma_{g,r}^n[l](x)$. We shall denote the subgroup
of $\Gamma_{g,r}^n(x)$ corresponding to $T_{g,r}^n$ by $T_{g,r}^n(x)$.
Denote the relative Malcev completion of $\Gamma_{g,r}^n(x)$ by
$\G_{g,r}^n(x)$, its prounipotent radical by $\U_{g,r}^n(x)$, etc.
The Lie algebras corresponding to $T_{g,r}^n(x)$ and $U_{g,r}^n(x)$
will be denoted by $\t_{g,r}^n(x)$ and $\u_{g,r}^n(x)$, respectively.

In this section we prove that for each choice of a point $x$ in
$\M_{g,r}^n$, there is a canonical $\Q$ mixed Hodge structure
(MHS) on $\t_{g,r}^n(x)$. The first ingredient in the construction
of this MHS is the following theorem, which is proved in
\cite[(13.1)]{hain:derham}.

\begin{theorem}\label{mhs_gen}
Suppose that $X$ is a smooth quasi-projective algebraic variety
and $(\V,\langle\blank,\blank\rangle)$ is a polarized variation
of Hodge structure over $X$ of geometric origin whose monodromy
representation
$$
\rho : \pi_1(X,x_o) \to \Aut_\R(V_o,\langle\blank,\blank\rangle)
$$
has Zariski dense image. Then the coordinate ring of the
completion of $\pi_1(X,x_o)$ relative to $\rho$ and its unipotent
radical both have natural real MHSs such that the product, coproduct,
and antipode of each are morphisms of MHSs.
\qed
\end{theorem}

We will say that a homomorphism $\G \to\H$ between proalgebraic groups,
each of whose coordinate rings is a Hopf algebra in the category of
mixed Hodge structures, is a morphism of MHSs if the corresponding map
on coordinate rings is.
Since $\Gamma_{g,r}^n(x)$ is the orbifold fundamental group of
$(\M_{g,r}^n,x)$, the following result is not unexpected.

\begin{theorem}
For all $g,r,n \ge 0$, and for each choice of a point
$$
x = [C;x_1,\dots,x_n;v_1,\dots,v_r]
$$
of $\M_{g,r}^n$, there is a canonical real MHS on the coordinate
ring of $\G_{g,r}^n(x)$ for which
the product, coproduct and antipode are morphisms of MHS. Moreover,
the homomorphisms $\G_{g,r}^n(x) \to \G_{g,r}^{n-1}(x')$ and
$\G_{g,r}^n(x) \to \G_{g,r-1}^{n+1}(x'')$, induced by forgetting a point or
by replacing a tangent vector by its anchor point, are morphisms
of mixed Hodge structure.
\end{theorem}

\begin{proof}
Since the mapping class group is not, in general, the fundamental group
of $\M_{g,r}^n$, we need to pass to a level.
Choose an integer $l$ such that $\Gamma_{g,r}^n[l]$ is torsion free. In
this case, the moduli space $\M_{g,r}^n[l]$ is smooth and has
fundamental group isomorphic to $\Gamma_{g,r}^n[l]$. Since
$\Gamma_{g,r}^n[l]$ is torsion free, there is a universal curve
$$
\pi : \cC \to \M_{g,r}^n[l].
$$
Take $\V$ to be the dual of the local system $R^1\pi_\ast \Z$. This is
a polarized variation of Hodge structure of weight $-1$ and is clearly
of geometric origin. Its monodromy representation is
$$
\rho : \Gamma_{g,r}^n[l] \to Sp_g(\Z).
$$
So by (\ref{level}) and (\ref{mhs_gen}), there is a canonical real
MHS on the coordinate ring of the relative completion $\G_{g,r}^n(x)$
of $\Gamma_{g,r}^n(x)$ for each choice of a point of $\M_{g,r}^n[l]$
that lies over $x$.

Denote the projection $\M_{g,r}^n[l] \to \M_{g,r}^n$ by $p$. The
set of lifts of a point $x$ of $\M_{g,r}^n$ to $\M_{g,r}^n[l]$ is
permuted transitively by the Galois group $Sp_g(\Z/l\Z)$. It follows
from the naturality of the MHS on the relative completion that the
MHSs on $\G_{g,r}^n(x)$ with respect any two points of $p^{-1}(x)$ are
canonically isomorphic. The MHS on $\Gamma_{g,r}^n(x)$ is therefore
indepenent of the choice of a point of $p^{-1}(x)$, and is therefore
canonical.

To show that the MHS on $\Gamma_{g,r}^n(x)$ constructed above is
independent of the choice of the level $l$, suppose that $l_1$ and
$l_2$ are two levels for which the mapping class group is torsion free.
One can then compare the corresponding MHSs by passing to the level
corresponding to the least common multiple of $l_1$ and $l_2$.
The naturality statement follows directly from
\cite[(13.12)]{hain:derham}.
\end{proof}

\begin{corollary}
For all $g\ge 0$, and for each choice of a point
$$
x = [C;x_1,\dots,x_n;v_1,\dots,v_r]
$$
of $\M_{g,r}^n$, the pronilpotent Lie algebra $\u_{g,r}^n(x)$ of the
prounipotent radical of $\G_{g,r}^n(x)$ has a canonical real MHS for
which the bracket is a morphism of MHS. Moreover, the morphisms
$\u_{g,r}^n(x) \to \u_{g,r}^{n-1}(x')$ and
$\u_{g,r}^n(x) \to \u_{g,r-1}^{n+1}(x")$,
obtained by forgetting a point or replacing a tangent vector by its
anchor point, are morphisms of mixed Hodge structure. \qed
\end{corollary}

Given a point $x$ of $\M_{g,r}^n$, there are two {\it a priori} different
MHSs on $\p_{g,r}^n(x)$. The first is the one obtained from the
construction given in \cite{hain:dht}. The second arises as
$\p_{g,r}^n(x)$ is the kernel of the natural surjection
$\u_{g,r}^n(x) \to \u_g(x)$. The following assertion follows directly
from the naturality properties \cite[(13.12)]{hain:derham} of the mixed
Hodge structure relative completions.

\begin{proposition}
These two MHSs are identical.\qed
\end{proposition}

Fix a point $x$ of $\M_{g,r}^n$. Then both of $\G_{g,r}^n(x)$ and
$\p_{g,r}^n(x)$ have canonical MHSs. It is natural to expect that the
natural action
\begin{equation}\label{action}
\G_{g,r}^n(x) \to \Aut \p_{g,r}^n(x)
\end{equation}
constructed in (\ref{rep}) is compatible with these.

\begin{lemma}\label{action_morph}
The action (\ref{action}) is a morphism of MHS. Consequently, the
morphism
$$
\u_{g,r}^n(x) \to \Der \p_{g,r}^n(x)
$$
is also morphism of MHS with respect to the canonical MHSs determined
by $x\in \M_{g,r}^n$.
\end{lemma}

\begin{proof}
It follows immediately from (\ref{exactness}) that $\cP_{g,r}^n(x)$
is a normal subgroup of $\G_{g,r}^n(x)$. Since the coordinate ring of
$\G_{g,r}^n(x)$ has a natural mixed Hodge structure compatible with its
operations, the action of $\G_{g,r}^n(x)$ on $\cP_{g,r}^n(x)$ via
conjugation is a morphism of MHS. But this action is easily seen to
coincide with the canonical action of $\G_{g,r}^n(x)$ on $\cP_{g,r}^n(x)$.
\end{proof}

For a curve $C$ of genus $g\ge 3$, denote by $PH_3(\Jac C,\Q)$ the
{\it primitive three dimensional homology} of its jacobian $\Jac C$
--- that is, the subspace of $H_3(\Jac C,\Q)$ corresponding to
$PH^{2g-3}(\Jac C,\Q)$ under Poincar\'e duality. It has a natural
Hodge structure of weight $-3$.

\begin{proposition}\label{purity}
If $g \ge 3$, then for each
$x = [C;x_1,\dots,x_n;v_1,\dots,v_r] \in \M_{g,r}^n$
the canonical real mixed Hodge structure on $H_1(\u_{g,r}^n(x))$ is of
weight $-1$ and is canonically isomorphic to
$$
PH_3(\Jac C,\R(-1)) \oplus H_1(C,\R)^{\oplus(r+n)}.
$$
\end{proposition}

\begin{proof}
As in the proof of (\ref{h1_tor}), we reduce to proof to showing
that it is true for $\u_g^1$. Then, by (\ref{action_morph}), the
composite
\begin{equation}\label{comp}
H_1(\u_g^1(x)) \to W_{-1}H_1(\Der \p_g^1(x)) \to \Gr^W_{-1}\Der\p_g^1(x)
\end{equation}
is a morphism of MHS. Observe that
$$
\Gr^W_{-1}\Der\p_g^1(x) \subset
\Hom(\Gr^W_{-1}\p_g^1(x),\Gr^W_{-2}\p_g^1(x)).
$$
{}From the work of Johnson \cite{johnson:def} (see also
\cite[\S4]{hain:normal}), it follows that (\ref{comp})
is injective, from which the result follows for $\u_g^1$.
\end{proof}

The fact that $H_1(\u_{g,r}^n)$ is pure of weight $-1$ allows us to
conclude that the weight filtration of $\u_{g,r}^n$ is essentially
its lower central series. This follows from the following general
fact.

\begin{lemma}\label{wt=lcs}
Suppose that $\g$ is a pronilpotent Lie algebra in the category of
mixed Hodge structures with finite dimensional $H_1$. If the induced
MHS on $H_1(\g)$ is pure of weight $-1$, then $W_{-l}\g$ is the $l$th
term of the lower central series of $\g$.
\end{lemma}

\begin{proof}
Denote the $l$th term of the lower central series of $\g$ by $\g^{(l)}$.
Since $H_1(\g)$ is pure of weight $-1$, since $\g$ is pronilpotent, and
since the bracket is a morphism of MHS, it follows that $\g = W_{-1}\g$.
An elementary argument using the Jacobi identity shows that the bracket
\begin{equation}\label{bra}
\g \otimes \g^{(l)} \to \g^{(l+1)}
\end{equation}
is surjective. Since the bracket is a morphism of MHS, it follows that
$\g^{(l)} \subseteq W_{-l}\g$. The fact that $H_1(\g)$ is pure of weight
$-1$, forces $\g^{(2)} = W_{-2}\,\g$. The result now follows by an induction
argument (induct on $l$) using the fact that (\ref{bra}), being a morphism
of MHS, is strict with respect to the weight filtration.
\end{proof}

\begin{corollary}\label{lcs}
The $l$th term of the lower central series of $\u_{g,r}^n$
is $W_{-l}\u_{g,r}^n$. \qed
\end{corollary}

This result implies that the weight filtration on $\G_{g,r}^n$
is defined over $\Q$ which implies that this MHS is really defined
over $\Q$.

\begin{corollary}
The weight filtration of the canonical MHSs on $\G_{g,r}^n(x)$ and
$\u_{g,r}^n(x)$ associated to  a point of $\M_{g,r}^n$ are topologically
determined and therefore defined over $\Q$. Consequently, the MHSs on
$\G_{g,r}^n(x)$ and $\u_{g,r}^n(x)$ each have a canonical lift to
$\Q$-MHSs.
\qed
\end{corollary}

We are now ready to lift the MHS from $\u_{g,r}^n(x)$ to $\t_{g,r}^n(x)$:

\begin{theorem}\label{mhs_torelli}
Suppose that $g\neq 2$ and that $r,n \ge 0$. For each choice of
a base point
$$
x = \left[C;x_1,\dots,x_n;v_1,\dots,v_r\right]
$$
of $\M_{g,r}^n$ there is a canonical $\Q$-MHS on $\t_{g,r}^n(x)$
for which the bracket and the quotient map $\t_{g,r}^n(x) \to u_{g,r}^n(x)$
are morphisms of MHS. Moreover, $W_{-l}\t_{g,r}^n(x)$ is the $l$th term
of the lower central series of $\t_{g,r}^n(x)$ and the central $\Ga$ is
isomorphic to $\Q(1)$.
\end{theorem}

\begin{proof}
For all $g \ge 0$, we have the exact sequence
$$
1 \to \pi_{g,r}^n \to T_{g,r}^n \to T_g \to 1.
$$
When $g =0,1$, $T_g$ is the trivial group, so that $T_{g,r}^n$ is
isomorphic to $\pi_{g,r}^n$. It follows that in these cases
$\t_{g,r}^n$ is  isomorphic to the Malcev Lie algebra $\p_{g,r}^n$
associated to $\pi_{g,r}^n$. The choice of the base point of
$\M_{g,r}^n$ gives the configuration space $(F_{g,r}^n,f_o)$ the
structure of a pointed smooth complex algebraic variety. Since
$\pi_{g,r}^n$ is the fundamental group of $(F_{g,r}^n,f_o)$, the
existence of the MHS on $\t_{g,r}^n(x)$ when $g=0,1$ follows from
\cite[(6.3.1)]{hain:dht}.

Now suppose that $g\ge 3$.
To construct a MHS on $\t_{g,r}^n(x)$, it suffices to show that $\t_g(x)$
has a MHS such that $\t_g(x) \to \u_g(x)$ is a morphism as it follows from
(\ref{central_ext}), (\ref{seq}) and (\ref{exactness}) that the diagram
$$
\begin{CD}
\t_{g,r}^n(x) @>>> \t_g(x) \cr
@VVV @VVV \cr
\u_{g,r}^n(x) @>>> \u_g(x) \cr
\end{CD}
$$
is a pullback square in the category of pronilpotent Lie algebras.

It is useful to begin by explaining the philosophy behind the proof.
The essential point is that $\Gamma_g(x)$ acts on $\t_g(x)$ and on
$\u_g(x)$ --- the action is induced by the action of $\Gamma_g(x)$
on $T_g(x)$ by conjugation. The central extension
\begin{equation}\label{extn}
0 \to \Ga \to \t_g(x) \to \u_g(x) \to 0.
\end{equation}
given by (\ref{central_ext}) can be viewed as an extension of
local systems over $\M_g$ (in the orbifold sense, of course)
where $\Ga$ is a trivial local system. Although we have not proved
it yet, $\u_g(x)$ should be a variation of MHS over $\M_g$. So we should
try to construct the MHS on $\t_g(x)$ so that (\ref{extn}) is both an
extension of local systems and an extension of mixed Hodge structures.
This, and the fact that the bracket has to be a morphism of MHS, gives
us no choice. We now carry this out this program.

The first point is that we know that, since $H_1(\t_g)$
does not contain any copies of the trivial representation,
the central $\Ga$ is contained in $[\t_g,\t_g]$.
The second is that by the computations in \S8 of
\cite{hain:comp}
we know that the central $\Ga$ lies in the image of the map
$$
\Lambda^2 H_1(\t_g) \to \Gr^{\mathrm{lcs}}_{-2} \t_g
$$
induced by the bracket. So in order that the bracket be a morphism
of MHS, the central $\Ga$ must be of weight $-2$. Since it is one
dimensional, it has to be isomorphic to $\Q(1)$.

Now fix a base point $x$ of $\M_g$. There is a corresponding $\Q$-MHS
on $\u_g(x)$. To lift this MHS to $\t_g(x)$, we have to give an element of
$$
\Ext^1_\H(\u_g(x),\Q(1)),
$$
where $\Ext_\H$ denote the Ext group in the category $\H$ of $\Q$
mixed Hodge structures. Applying the functor $\Ext_\H^\dot$ to the
sequence
$$
0 \to W_{-2} \u_g(x) \to \u_g(x) \to H_1(\u_g(x)) \to 0,
$$
we see that the natural map
$$
\Ext^1_\H(H_1(\u_g(x)),\Q(1)) \to \Ext^1_\H(\u_g(x),\Q(1))
$$
is an isomorphism.

Now let the base point vary. By (\ref{purity}), $H_1(\u_g)$
is a variation of $\Q$-MHS over $\M_g$ (in the orbifold sense)
of weight $-1$ --- cf.\ \cite[(9.1)]{hain:normal}. Denote the category
of admissible variations of $\Q$-MHS over a smooth variety $X$ by
$\H(X)$ and the category of $\Q$ local systems over $X$ by $\cL(X)$.
Then, by \cite[(8.4)]{hain:normal}, the the forgetful map
\begin{equation}\label{ext_iso}
\Ext^1_{\H(\M_g)}(H_1(\u_g),\Q(1))
\to \Ext^1_{\cL(\M_g)}(H_1(\u_g),\Q(1))
\cong H^1(\Gamma_g,H^1(\u_g))
\end{equation}
is an isomorphism.\footnote{This group is one dimensional and
generated by the Johnson homomorphism, although we do not really
need to know this here --- see \cite[(5.2)]{hain:normal}.}
We will lift the MHS on $\u_g(x)$ to $\t_g(x)$ using an element of the
right hand group.\footnote{The class we seek, not surprisingly, is the
one corresponding to the Johnson homomorphism, and is half the class
associated to the algebraic cycle $C-C^-$ in $\Jac C$ --- see
\cite[\S8]{hain:normal}.} We do this by producing an element of the
right hand group which corresponds to the central extension (\ref{extn}).

Denote the $k$th term of the lower central series of $\t_g(x)$
by $\t_g(x)^{(k)}$. We can form the extension
\begin{equation}\label{extn2}
0 \to \t_g(x)^{(2)}/\t_g(x)^{(3)} \to
\t_g(x)/\t_g(x)^{(3)} \to  H_1(\t_g(x)) \to 0
\end{equation}
As has been pointed out above, the image of $\t_g(x)^{(2)}/\t_g(x)^{(3)}$
in the central $\Ga$ in $\T_{g,r}^n(x)$ is non-trivial. The kernel of
the extension (\ref{extn2}) is a rational
representation of $Sp_g$. Since $H_1(\t_g(x))$ is irreducible, its
second exterior power contains exactly one copy of the trivial
representation. There is therefore a unique non-zero $Sp_g$-invariant
projection
$$
\t_g(x)^{(2)}/\t_g(x)^{(3)} \to \Ga.
$$
If we push the extension (\ref{extn2}) out along this map,
we obtain a canonical element of
$$
\Ext^1_{\cL(\M_g)}(H_1(\u_g),\Q(1)).
$$
By the isomorphism (\ref{ext_iso}), we obtain an element
of
$$
\Ext^1_{\H(\M_g)}(H_1(\u_g),\Q(1))
$$
which allows us, for each point $x$ of $\M_g$, to lift the MHS
on $\u_g(x)$ associated to $x$ to a MHS on $\t_g(x)$.

Our last task is to show that the bracket is a morphism of MHS.
We only need show that the bracket preserves the Hodge and
weight filtrations.
First observe that since $\Ga$ is central and contained in $W_{-2}$,
the bracket preserves the weight filtration, and its restriction to
$W_{-2}\t_g(x)$ is a morphism of MHS. It remains to show that the
the bracket preserves the Hodge filtration. In view of these facts,
it suffices to prove that
$$
\left[F^p\t_g(x),F^q\t_g(x)\right] = 0
$$
when $p+q > -1$. This is easily deduced from the the
fact that the map
$$
[\blank,\blank] : \Lambda^2 H_1(\t_g)
\to \Gr^W_{-2} \t_g(x) \stackrel{\text{proj}}{\to} \Ga\cong \Q(1)
$$
induced by the bracket is a polarization of $H_1(\t_g(x))$ as it is
$Sp_g$ equivariant and non-zero by results in \cite[\S7]{hain:comp}.
\end{proof}

\begin{remark}
It follows immediately from (\ref{action_morph}) and
(\ref{mhs_torelli}) that for  each choice of base point in
$\M_{g,r}^n$, the canonical morphism
$\t_{g,r}^n(x) \to \Der \p_{g,r}^n(x)$ is a morphism of MHS.
\end{remark}

\section{Review of Continuous Cohomology}
\label{cts_coho}

In this section, we briefly review the theory of continuous cohomology
of discrete groups, which is mainly developed in \cite{hain:cycles}.
It will be our principal tool in proving that Torelli has a
presentation with only quadratic relations. As a warm up, we show
how it can be used to give a new and simpler proof of Morgan's theorem
that the complex form of the Lie algebra associated to the fundamental
group of a smooth variety has a weighted homogenous presentation with
generators of weights equal to those occurring in $H_1(X)$, and relations
of weight contained in those of $H_2(X)$.

Define the continuous cohomology of a discrete group $\pi$ to be
the direct limit of the rational cohomology of its finitely
generated nilpotent quotients:
$$
\Hcts^\dot(\pi,\Q) := \lim_\to H^\dot(N,\Q) \label{ctsgp_def}
$$
where $N$ ranges over the finitely generated nilpotent quotients
of $\pi$. There is an obvious natural homomorphism
\begin{equation}\label{natural_homom}
\Hcts^\dot(\pi,\Q) \to H^\dot(\pi,\Q).
\end{equation}
If $X$ is a topological space with fundamental group $\pi$, then
we also have a natural homomorphism
$$
\Hcts^\dot(\pi,\Q) \to H^\dot(X,\Q)
$$
as there is a canonical map $H^\dot(\pi) \to H^\dot(X)$.

\begin{proposition}\label{cts_ord}
If $H_1(\pi,\Q)$ is finite dimensional, then the natural homomorphism
(\ref{natural_homom}) is an isomorphism in degree 1 and injective in
degree 2. \qed
\end{proposition}

This is really a restatement of the result of Dennis Sullivan which
asserts that the Lie algebra of the 1-minimal model of a space is
the Malcev Lie algebra of the fundamental group. A more direct proof
can be found, for example, in \cite[(5.1)]{hain:cycles}. We present a
new proof because it is elementary.

\begin{proof}
The group $H^2(G,\Q)$ parameterizes central extensions of $G$ by $\Q$.
Suppose that $\alpha \in \Hcts^2(\pi,\Q)$ whose image in $H^2(\pi,\Q)$
is trivial. Then there is a nilpotent quotient $N$ of $\pi$ and an element
$\tilde{\alpha}$ of $H^2(N,\Q)$ that is a lift of $\alpha$. There is a
central extension
\begin{equation}\label{alpha}
1 \to \Q \to E \to N \to 1
\end{equation}
corresponding to  by $\tilde{\alpha}$. The key point to note is that $E$ is
nilpotent. To say that the image of $\alpha$ is trivial in $H^2(\pi,\Q)$
is to say that the pullback of the extension (\ref{alpha}) to $\pi$ is
split. Composing a splitting of this projection with the projection of
$\pi \to N$ gives  a homomorphism $\pi \to E$ which lifts $\pi \to N$.
Denote the image of $\pi$ in $E$ by $\tilde{N}$. It is easy to see that
the pullback of the extension (\ref{alpha}) to $\tilde{N}$ splits. Since
$\tilde{N}$ is a nilpotent quotient of $\pi$, the class $\alpha$ vanishes.
\end{proof}

Similarly, we can define the continuous cohomology of a pronilpotent
Lie algebra $\g$ to be the direct limit of the
cohomology of its finite dimensional nilpotent quotients:
$$
\Hcts^\dot(\g) := \lim_\to H^\dot(\n) \label{ctslie_def}
$$
where $\n$ ranges over the finite dimensional nilpotent quotients
of $\g$.

A mild generalization of a theorem of Nomizu \cite{nomizu} states
that for each finitely generated nilpotent group $N$ there is a natural
isomorphism
$$
H^\dot(\n) \cong H^\dot(N,\Q)
$$
where $\n$ is the Lie algebra of the Malcev completion of $N$. It
follows immediately from the definitions that if $\pi$ is a finitely
generated group and $\p$ the associated Malcev Lie algebra, then there
is a natural isomorphism
$$
\Hcts^\dot(\p) \cong \Hcts^\dot(\pi,\Q).
$$

The continuous cohomology of a pronilpotent Lie algebra $\g$ can
be computed using the standard complex $\cC^\dot(\g)$ of continuous
cochains of $\g$. This is defined to be the direct limit of
the Chevalley-Eilenberg cochains of the finite dimensional nilpotent
quotients of $\g$. Denote the continuous dual of $\g$ by $\g^\ast$.
Then we have a d.g.a.\ isomorphism
$$
\cC^\dot(\g) = \Lambda^\dot \left(\g^\ast[-1]\right);
$$
the differential is derivation of degree 1 whose restriction
to $\g^\ast$ is minus the dual of the bracket.

The definition of continuous cohomology can be extended to the case
where the coefficients are $\Q$ modules on which $\g$ acts via
a representation of one of its nilpotent quotients --- cf.\
\cite{hain:cycles}.

Suppose now that $H_1(\g)$ is finite dimensional. If $\g$ has a MHS,
then, by linear algebra, so do $\cC^\dot(\g)$ and $H^\dot(\g)$. We will
call such a pronilpotent Lie algebra a {\it Hodge Lie algebra}. It
follows that if $X$ is an algebraic variety, $x\in X$, then
$\Hcts^\dot(\pi_1(X,x),\Q)$ has a canonical MHS. One can show that
this MHS does not depend on the base point $x$ of $X$ ---
\cite{hain:cycles}.

Since the weight filtration of a MHS splits canonically over $\C$,
each finite dimensional Hodge Lie algebra $\g$ is canonically
isomorphic to the graded Lie algebra $\Gr^W_\dot\g$ after tensoring
with $\C$. The following result therefore follows by taking inverse
limits.

\begin{proposition}\label{canon_split}
If $\g$ is a Hodge Lie algebra, all of whose weights are negative,
then there is a canonical Lie algebra isomorphism
$$
\g_\C \cong \prod_{l\ge 1} \Gr^W_{-l}\g_\C. \qed
$$
\end{proposition}

Since each choice of a base point of $\M_{g,r}^n$ determines a
canonical MHS on $\t_{g,r}^n$, we have:

\begin{corollary}
For each choice of a base point of $\M_{g,r}^n$, there is a
canonical isomorphism
$$
\t_{g,r}^n\otimes \C \cong
\prod_{l\ge 1}\left(\Gr^W_{-l} \t_{g,r}^n \otimes \C\right). \qed
$$
\end{corollary}

The following result is proved, for example, in
\cite[(11.7)]{carlson-hain}. In Section \ref{cts_coho_tor} we will
prove a generalization needed for studying the relations in $\t_g$.

\begin{theorem}\label{cts_morph}
If $X$ is a smooth algebraic variety, then the natural homomorphism
$$
\Hcts^\dot(\pi_1(X),\Q) \to H^\dot(X,\Q)
$$
is a morphism of mixed Hodge structures.
\end{theorem}

The final two results in this section together will allow us to
use continuous cohomology as an effective tool for studying relations
in Hodge Lie algebras in general, and $\t_g$ in particular.

The cochains, and therefore the cohomology, of a graded Lie algebra
both have an extra grading, and are therefore bigraded algebras. If
$\g$ is a Hodge Lie algebra, then $\Gr^W_\dot\g$ has an extra grading
is by weight. Since the functor $\Gr^W$ is exact on the category of
MHS, we have:

\begin{proposition}
If $\g$ is a Hodge Lie algebra, then there is a
canonical bigraded algebra isomorphism
$$
\Gr^W_\dot \Hcts^\dot(\g) \cong H^\dot(\Gr^W_\dot \g). \qed
$$
\end{proposition}

If $\g$ is a graded Lie algebra with negative weights, then we can write
$\g$ as a quotient of the free graded Lie algebra $\f$ generated by
$H_1(\g)$ modulo a homogeneous ideal $\r$. Note that we are not assuming
that $H_1(\g)$ is pure --- in general it will be graded. The group
$$
H_0(\f/\r) = \r/[\f,\r]
$$
is graded. One can obtain a minimal set of relations of $\g$ by
taking the image of any splitting of the projection
$$
\r \to H_0(\f/\r).
$$
The following result is an analogue of Hopf's description
of the second homology of a group in terms of a presentation.

\begin{proposition}
If $\g$ is a graded Lie algebra with negative weights,
then there is a canonical isomorphism of graded vector spaces
$$
H_0(\f/\r) \cong H_2(\g).
$$
\end{proposition}

\begin{proof}
There are  several ways to see this. One is look closely at the
Chevalley-Eilenberg cochains of $\g$. The second is the use the
fact that a sub-Lie algebra of a free Lie algebra is free \cite%
[(2.5)]{reutenauer} to deduce that, as a Lie algebra, $\r$ is free.
Then apply the Lie algebra analogue of the Hochschild-Serre spectral
sequence to the extension
$$
0 \to \r \to \f \to \g \to 0.
$$
The details are standard and are omitted.
\end{proof}

\begin{corollary}\label{gr_presentn}
If $\g$ is a graded Lie algebra with negative weights, then
there is an injective linear map
$$
\delta : H_2(\g) \hookrightarrow \L(H_1(\g))
$$
of graded vector spaces such that $\g$ has presentation
$$
\L(H_1(\g))/(\im \delta)
$$
in the category of graded Lie algebras.
\end{corollary}

Combining (\ref{canon_split}), (\ref{gr_presentn}) and the
existence of a canonical MHS on the Malcev Lie algebra $\g(X,x)$
associated to a pointed variety, we obtain
one of Morgan's theorems \cite[(10.3)]{morgan}.

\begin{theorem}\label{morgan}
If $X$ is a smooth complex algebraic variety and $x\in X$, then
then the complex Malcev Lie algebra $\g(X,x)_\C$ associated to
$\pi_1(X,x)$ has the property that
$$
\g(X,x)_\C \cong \prod_{l \ge 1}\Gr^W_{-l}\g_C
$$
and there is a homomorphism of graded vector spaces
$$
\delta : H_2(X,\C) \to \L(\Gr^W_\dot H_1(X))
$$
such that
$$
\Gr^W_\dot \g_\C \cong
\L(\Gr^W_\dot H_1(X,\C))/(\delta(\Gr^W_\dot H_2(X,\C)))
$$
in the category of graded Lie algebras. \qed
\end{theorem}

\section{Remarks on the Representations of $\sp_g$}
\label{reps}

In this section we review some basic facts from the representation
theory that we shall need in subsequent sections. A
good reference is \cite{fulton-harris}.

Denote the Lie algebra of $Sp_g$ by \label{symp_def} $\sp_g$. The
representation theory of the group and the Lie algebra are the same.
Denote their common representation ring by \label{rep_def} $R(\sp_g)$.

Choose a symplectic basis $a_1,\dots,a_g,b_1,\dots,b_g$ of the
fundamental representation of $\sp_g$. Denote by $\h$ the torus in
$\sp_g$ consisting of matrices that are diagonal with respect to this
basis. Choose coordinates $t=(t_1,\ldots,t_g)$ on $\h$ so that
$$
t\cdot a_i = t_ia_i \text{ and } t\cdot b_i = - t_i b_i.
$$
The subalgebra of positive nilpotents $\n$ has basis the elements
$S_{i,j}, (i<j)$, $T_i$, and $F_{i,j}, (i\neq j)$ of $\sp_g$,
where
$$
S_{i,j} (a_j) = a_i, \quad S_{i,j}(b_i) = - b_j,  \quad
S_{i,j}(\text{other basis vectors}) = 0,
$$
$$
T_i(b_i) = a_i, \quad T_i(\text{other basis vectors}) = 0,
$$
$$
F_{i,j}(b_i) = a_j, \quad F_{i,j} (b_j) = a_i,  \quad
F_{i,j}(\text{other basis vectors}) = 0.
$$

A fundamental set of weights of $\sp_g$ is $\lambda_j:\h \to \R$,
$1\le j \le g$, where $\lambda_j$ is defined by
$$
\lambda_j(t) = t_1 + t_2 + \cdots + t_j. \label{wt_def}
$$
The irreducible representations of $\sp_g$ correspond to
positive integral linear combinations $\lambda$ of the $\lambda_j$.
Denote the irreducible representation of $\sp_g$ with highest
weight $\lambda$ by \label{module_def} $V(\lambda)$.

The irreducible representations of $\sp_g$ can also be indexed
by partitions $\alpha$ of an integer $n$ into $\le g$ parts:
$$
n = \alpha_1 + \alpha_2 + \dots + \alpha_g
$$
where
$$
\alpha_1 \ge \alpha_2 \ge \dots \ge \alpha_g \ge 0.
$$
The irreducible representation corresponding to $\alpha$
has highest weight
$$
t \mapsto \sum_j \alpha_j t_j.
$$

We shall denote the integer
$$
\sum_{k=1}^g \alpha_k = \sum_{k=1}^g k\, n_k \label{size_def}
$$
by $|\alpha|$ or by $|\lambda|$ according to whether we are using
partitions or highest weights. This can be considered as a measure
of the size of the corresponding irreducible representation; it is the
smallest positive integer $d$ such that $V(\lambda)$ occurs in the
$d$th tensor power of the fundamental representation.

There is a notion of stability of the decomposition of tensor
products and Schur functors of representations of symplectic groups.
In order to state the result, we need to first define the {\it depth},
$\delta(V)$, of a representation $V$ of $\sp_g$. If the module is
irreducible with highest weight $\sum n_k\,\lambda_k$, define
$\delta(V)$ to be the largest $d$ such that $n_d\neq 0$ --- or
equivalently, it is the number of rows in the corresponding Young
diagram. Define the depth of an arbitrary representation to be the
maximum of the depths of its irreducible components.

In order to discuss stability, we will need a stabilization map.
When $h \ge g$, define a group homomorphism
$$
R(\sp_g) \hookrightarrow R(\sp_h)
$$
by taking the irreducible representation of $\sp_g$ corresponding to the
partition $\alpha$ to the representation of $\sp_h$ corresponding to the
same partition. Equivalently, take the representation of $\sp_g$ with
highest weight $\sum n_k\, \lambda_k$ to the representation of $\sp_h$
with the same highest weight decomposition.

Recall that to each partition $\beta$ of a positive integer $n$,
one has a Schur functor $\Schur_\beta$ defined on the category of
representations of each group. For example, if $\beta = [n]$,
then $\Schur_\beta$ is the $n$th symmetric power, if $\beta = [1^n]$,
then $\Schur_\beta$ is the $n$th exterior power. We shall denote the
integer $n$ by $|\beta|$.

The second assertion of the following stability result appears to be
folklore --- the only proof I know of is in Kabanov's thesis.

\begin{theorem}\label{kab_stab}
\begin{enumerate}
\item (\cite[p.~424]{fulton-harris}) If $V$ and $W$ are representations
of $\sp_g$ and $\delta(V)
+ \delta(W) \le g$, then the irreducible representations and their
multiplicities occurring in the decomposition of $V\otimes W$ is
independent of $g$.
\item (\cite{kabanov} --- see also \cite{kabanov:stab}) If $V$ is a
representation of $\sp_g$ and $\beta$ is a partition
with $|\beta|\delta(V) \le g$, then the decomposition of
$\Schur_\beta V$ into irreducible components is independent of $g$.
\end{enumerate} \qed
\end{theorem}

\begin{remark}\label{method_comp}
Some of the computations of highest weight decompositions in this
paper have been made using the computer program \LiE\ from the
University of Amsterdam. The computations were performed for a
particular $g$ in the stable range. The stability theorem was then
used to deduce the decomposition for all $g$ in the stable range. Note
that all such computations were checked using several values of $g$ in
the  stable range. In addition, many of the unstable computations were
done using \LiE.
\end{remark}

By composition with the canonical homomorphism
$$
\Gamma_{g,r}^n \to Sp_g(\Q)
$$
we see that each representation $V$ of $\sp_g$ gives rise to a local
system over $\M_{g,r}^n$, at least in the orbifold sense.
It is a standard fact that each such local system arising from an
irreducible representation of $\sp_g$ is an admissible variation
of Hodge structure over $\M_{g,r}^n$ in a unique way up to Tate
twist  --- cf.\ \cite[(9.1)]{hain:normal}. It can always be realized
as a variation of weight $|\lambda|$, and we shall take this as the
default weight.

We would like to discuss the cohomology of $\M_{g,r}^n$ with
coefficients in such a local system. To do this, first choose
a level $l$ such that $\Gamma_{g,r}^n$ is torsion free. Then
$\M_{g,r}^n[l]$ is smooth and the variation of Hodge structure $\V$
corresponding to an irreducible representation $V$ of $Sp_g$ is defined
over $\M_{g,r}^n[l]$ and has a natural $Sp_g(\Z/l\Z)$ action. From the
work of M.~Saito \cite{saito}, we know that $H^k(\M_{g,r}^n[l],\V)$ has
a canonical mixed Hodge structure with weights
$\ge k + \text{weight}(V)$. The action of $Sp_g(\Z/l\Z)$ preserves
this MHS. So we can define
$$
H^\dot(\M_{g,r}^n,\V) = H^0(Sp_g(\Z/l\Z),H^\dot(\M_{g,r}^n[l],\V)).
$$
as a MHS. Note that the underlying group is canonically isomorphic to
$H^\dot(\Gamma_{g,r}^n,V)$.

\section{Continuous Cohomology of Torelli Groups}
\label{cts_coho_tor}

The next step in finding a presentation of $\t_{g,r}^n$ is to determine
the relations in $\Gr^W_\dot \t_{g,r}^n$. Since this is a graded Lie
algebra generated in degree $-1$, the generators of the ideal of
relations is homogeneous. In this section we will use a result of
Kabanov \cite{kabanov} (see also \cite{kabanov:purity}) about the second
cohomology of $\Gamma_{g,r}^n$
to show that the ideal of relations in $\Gr\t_{g,r}^n$ is generated
by quadratic and cubic generators when $g\ge 3$, and quadratic
relations alone when $g\ge 6$.
Our principal tool will be the continuous cohomology defined in
Section \ref{cts_coho}.

First some notation. Take $X$ and $\V$ as in the
statement of (\ref{mhs_gen}). Denote the Lie algebra associated to
the prounipotent radical of the completion of $\pi_1(X,x)$ relative
to $\rho$ by $\u(x)$. This is a Hodge Lie algebra. The next result is
a generalization of (\ref{cts_morph}).

\begin{proposition}\label{morphism}
Suppose that $\W$ is an admissible variation of Hodge structure over
$X$ which is a subquotient of a tensor power of $\V$.  Then for all
$k\ge 0$ and each $x\in X$, there is a natural homomorphism
$$
H^0(X,\Hcts^k(\u(x))\otimes W_x) \to H^k(X,\W)
$$
which is a morphism of MHS. It is an isomorphism when $k=1$ and
injective when $k=2$.
\end{proposition}

\begin{proof}
The case $k=1$ is proved in \cite[(10.3),(13.8)]{hain:derham}.
We will prove the result when $k>1$ by induction. The most important
case for us is when $k=2$, so we will give that argument
in more detail and briefly sketch the remaining cases. We will assume
throughout that the reader is familiar with \cite{hain:derham}. A
convenient auxiliary reference for rational homotopy theory is
\cite[\S2]{hain:dht}.

We begin by recalling some well known facts from rational
homotopy theory. The base point $x\in X$ determines an
augmentation
$$
\epsilon_x : \Efin^\dot(X,\O(P)) \to \R,
$$
where $P$ is the principal bundle defined in \cite[\S 4]{hain:derham}.
We shall write $\O$ instead of $\O(P)$.
We can form the bar construction
$$
B(\Efin^\dot(X,\O)) := B(\R,\Efin^\dot(X,\O),\R)
$$
where both copies of $\R$ are regarded as $\Efin^\dot(X,\O)$
modules via $\epsilon_x$. The Lie algebra $\u(x)$ is determined
by $B(\Efin^\dot(X,\O))$ as follows: the dual
$$
H^0(B(\Efin^\dot(X,\O)))
$$
is a complete Hopf algebra, $\u(x)$ is its set of primitive elements.
(See, for example, \cite[(2.4.5) and \S 2.6]{hain:dht}.) There is an
augmentation preserving d.g.a.\ homomorphism
\begin{equation}\label{min_model}
\cC^\dot(\u(x)) \to \Efin^\dot(X,\O),
\end{equation}
unique up to homotopy, which induces the map
$$
\theta : \Hcts^\dot(\u(x)) \to H^\dot(\Efin^\dot(X,\O))
$$
on homology. The map $\theta$ is an isomorphism in degree 1 and injective
in degree 2.\footnote{In the language of Sullivan \cite{sullivan}, the
map (\ref{min_model}) is the 1-minimal model of $\Efin^\dot(X,\O)$.}
There is a canonical isomorphism
$$
H^k(X,\W) \cong H^0(X,H^k(\Efin^\dot(X,\O))\otimes\W)
$$
of MHSs for each VHS $\W$ over $X$ whose monodromy representation
is the pullback of a rational representation of $\Aut(V_o,q)$ via
the representation $\rho$. Since $\u(x)$ has a canonical MHS, and
since $\Efin^\dot(X,\O)$ is a mixed Hodge complex, each of the domain
and target of $\theta$ have a canonical MHS.  To prove
the result, we need only prove that $\theta$ is a morphism
of MHS.

First we give an intuitive proof. The image of the map
$$
\theta^2 : \Hcts^2(\u) \to H^2(\Efin^\dot(X,\O))
$$
is the subspace of the right hand side generated by the cup product
$H^1\otimes H^1 \to H^2$, all Massey triple products of 1-forms,
all Massey quadruple products of 1-forms, etc. Since the cup product
and all Massey $k$-fold products have domain which is a sub-MHS of
$\otimes^k H^1$ and are themselves
morphisms, it follows that the image of $\theta^2$ is a MHS. That
$\theta^2$ is a morphism follows as $\theta^1$ is an isomorphism of MHS.
One can  continue in an analogous fashion to prove that each $\theta^k$
is a morphism.

We now make this argument precise. The spectral sequence
associated to the standard filtration of the bar construction
is called the Eilenberg-Moore spectral sequence (EMss): for an
augmented d.g.a.\ $A^\dot$ with connected homology, it takes the
form
$$
E_1^{-s,t} = \left[\otimes^s H^+(A)\right]^t \implies
H^{t-s}(B(A)).
$$
Denote the EMss associated to $\cC(\u)^\dot$ by $\{E_r(\u)\}$ and the
EMss associated to $\Efin^\dot(X,\O)$ by $\{E_r(X)\}$.

The map (\ref{min_model}) induces a morphism of Eilenberg-Moore
spectral sequences. Each of these is a spectral sequence of MHSs
as both the domain and target of (\ref{min_model}) are mixed Hodge
complexes, but we have to prove that the map between them is a morphism
of MHSs. This is the case in total degree 0 as $E_1^{-s,s}$ is
$\otimes^s H^1$ and $\theta^1$ is an isomorphism of MHS.

It is a standard fact that
$$
H^k(B(\cC^\dot(\u))) = 0
$$
when $k>0$; cf.\ \cite[(2.6.2)]{hain:dht} and \cite{bloch-kriz}.
Therefore, the $E^{-1,2}_\infty$ term of the associated EMss vanishes.
(This is a precise way to say that $\Hcts^2(\u)$ is generated by Massey
products.) The edge homomorphisms
$$
\Hcts^k(\u) = E_1^{-1,k}(\u) \to E_r^{-1,k}(\u)
$$
are all surjective. Let $M_r^k$ be the inverse image in $\Hcts^k(\u)$
of the image of
$$
d_{r-1} : E_{r-1}^{-r,k+r-2}(\u) \to E_{r-1}^{-1,k}(\u).
$$
Then the fact that the higher cohomology of $B(\cC^\dot(\u))$ vanishes
implies that whenever $k\ge 2$
$$
\Hcts^k(\u) = \bigcup_r M_r^k.
$$
Since the spectral sequence is a spectral sequence of MHS, each
$M_r^k$ is a sub-MHS of $\Hcts^k(\u)$.

Since both spectral sequences are spectral sequences of MHSs, it
follows that the image of
$$
\Hcts^2(\u)=E_1^{-1,2}(\u) \to E_1^{-1,2}(X) = H^2(\Efin^\dot(X,\O))
$$
is a sub-MHS and that $\theta^2$ is a morphism of MHS.

If $k>2$, one can assume by induction that $\theta^m$ is a morphism
whenever $m<k$. It follows easily that the natural map
$$
E_1^{-s,t}(\u) \to E_1^{-s,t}(X)
$$
is a morphism of MHS whenever $-s+t < k-1$, and therefore that its
image is a sub-MHS of $E_1^{-s+t}(X)$. But since these spectral
sequences are spectral sequences in the category of MHSs, and since
$E_\infty^{-1,k}(\u)$ vanishes, it follows that $\theta^k$ is a
morphism.
\end{proof}

\begin{remark}
This is a continuation of Remark~\ref{sl2}. The previous result implies
that $\U_1[l]$ is a free pronilpotent group as $SL_2(\Z)[l]$ has a free
subgroup of finite index, which implies that $H^2(\M_1[l],S^n V)$
vanishes for all $n$. It follows that $H^2(\u_1[l])$ vanishes and from
(\ref{gr_presentn}) that $\u_1[l]$ is free.
\end{remark}

We thus have the following version of (\ref{morphism}) for moduli
spaces of curves.

\begin{proposition}\label{morph_u}
If $g\ge 3$ and $\V$ is a variation of Hodge structure over $\M_{g,r}^n$
whose monodromy representation comes from a rational representation
of $Sp_g$, then for all $k$, there is a natural map
$$
H^0(Sp_g,\Hcts^k(\u_{g,r}^n(x))\otimes V_x) \to
H^k(\M_{g,r}^n,\V)
$$
which is a morphism of MHS. Here $V_x$ denotes the fiber of $\V$
over $x$. \qed
\end{proposition}

This yields the following useful result about
differentials in the Hochschild-Serre spectral sequence associated
to the group extension
\begin{equation}\label{std}
1 \to T_{g,r}^n \to \Gamma_{g,r}^n \to Sp_g(\Z) \to 1.
\end{equation}
If we take coefficients in the irreducible representation $V(\lambda)$
of $Sp_g$, this spectral sequence takes the form
$$
E_2^{s,t} = H^0(Sp_g(\Z),H^t(T_{g,r}^n)\otimes V(\lambda))
\implies H^{s+t}(\Gamma_{g,r}^n,V(\lambda)).
$$

\begin{corollary}\label{vanishing}
For each $\lambda$, the image of the composite
$$
H^0(Sp_g,\Hcts^2(\t_{g,r}^n)\otimes V(\lambda))
\to H^0(Sp_g(\Z),H^2(T_{g,r}^n)\otimes V(\lambda)) = E_2^{0,2}
$$
is contained in $E_\infty^{0,2}$.
\end{corollary}

\begin{proof}
The result follows immediately from the fact that the diagram
$$
\begin{CD}
H^0(Sp_g,\Hcts^2(\u_{g,r}^n)\otimes V(\lambda)) @>>>
H^0(Sp_g,\Hcts^2(\t_{g,r}^n)\otimes V(\lambda)) \cr
@VVV @VVV \cr
H^2(\Gamma_{g,r}^n,V(\lambda)) @>>>
H^0(Sp_g(\Z),H^2(T_{g,r}^n)\otimes V(\lambda))
\end{CD}
$$
commutes, where the top map is the surjection induced by the projection
of $t_{g,r}^n$ onto $\u_{g,r}^n$, the right hand vertical map by the
canonical map
$$
\Hcts^\dot(\t_{g,r}^n) \to H^\dot(T_{g,r}^n)
$$
described in Section \ref{cts_coho}, the bottom map is the canonical
restriction map, and the left hand vertical map is the one given by
(\ref{morph_u}). This assertion can be proved using the constructions
in \cite[\S4]{hain:derham} by restricting to a leaf in the
principal bundle $P \to \M_{g,r}^n$ associated to the representation
(\ref{map}). In this case, each leaf is a copy of the classifying
space of $T_{g,r}^n$.
\end{proof}

Actually, we have proved a stronger statement than asserted. The
stronger claim will be stated in \S\ref{applications}.

Denote the fiber over $x\in \M_{g,r}^n$ of the variation of Hodge
structure $\V(\lambda)$ corresponding to the irreducible representation
$V(\lambda)$ of $Sp_g$ by $V(\lambda)_x$.

\begin{corollary}
If $g\ge 3$, then for each irreducible representation $V(\lambda)$
of $Sp_g$, there is a canonical monomorphism of MHS
$$
H^0(Sp_g,\Hcts^2(\t_{g,r}^n(x))\otimes V(\lambda)_x)
\hookrightarrow H^2(\M_{g,r}^n,\V(\lambda)).
$$
\end{corollary}

\begin{proof}

We first prove the existence of the homomorphism. Injectivity will
then follow directly from (\ref{morphism}). It follows from
(\ref{central_ext}) and \cite[(5.5)]{hain:cycles} that the sequence
of $Sp_g$ modules
$$
0 \to \Q(1) \to \Hcts^2(\u_{g,r}^n(x)) \to \Hcts^2(\t_{g,r}^n(x)) \to 0
$$
is an exact sequence of MHSs. It follows that the
natural map $\Hcts^2(\u_{g,r}^n(x)) \to \Hcts^2(\t_{g,r}^n(x))$ is a
surjective morphism of MHS. According to (\ref{morph_u}), the map
$$
\left[\Hcts^2(\u_{g,r}^n(x))\otimes V(\lambda)_x\right]^{Sp_g}
\to H^2(\M_{g,r}^n,\V(\lambda))
$$
is a morphism of MHS. So to construct the homomorphism, it suffices to
show that this map factors through the quotient map
$$
\Hcts^2(\u_{g,r}^n) \to \Hcts^2(\t_{g,r}^n).
$$

We first assume that $\lambda \neq 0,\lambda_1,\lambda_3$. Consider
the Hochschild-Serre spectral sequence
$$
E_2^{s,t} = H^s(Sp_g(\Z),H^t(T_{g,r}^n)\otimes V(\lambda))
\implies H^{s+t}(\Gamma_{g,r}^n,V(\lambda)).
$$
By (\ref{imp_borel}), $E_2^{s,t}$ vanishes when $s\le 1$ and
$t\le 2$ provided $g \ge 3$ (cf.\ \cite[(5.2)]{hain:normal}). It follows
that
$$
H^2(\Gamma_{g,r}^n,\V(\lambda)) =
H^0(Sp_g(\Z),H^2(T_{g,r}^n)\otimes V(\lambda)).
$$
The result now follows from (\ref{cts_ord}).

When $\lambda=\lambda_3$, we have $E_2^{2,1} \cong \Q$ (cf.\
\cite[(5.2)]{hain:normal}), so there is a possibility of having
a non-trivial differential $d_2 : E_2^{0,2} \to E_2^{2,1}$. But
by (\ref{vanishing}) this cannot occur. The argument is completed
as in the previous case. The case of $\lambda_1$ is proved in the
same way.

Finally, we consider the case of the trivial representation. In this
case, we have $E_2^{2,1}=E_2^{3,0}=0$, but $E_2^{2,0}=\Q$. It
follows that we have an exact sequence
$$
0 \to \Q \to H^2(\Gamma_{g,r}^n,\Q) \to
H^0(Sp_g(\Z),H^2(T_{g,r}^n)) \to 0.
$$
The result in this case now follows using the exact sequence in
the first paragraph of this proof.
\end{proof}

Denote the $\lambda$ isotypical part of an $Sp_g$ module $V$
by \label{iso_def} $V_\lambda$.

\begin{corollary} If $g\ge 3$ and $\lambda$ is a dominant integral
weight of $Sp_g$, then
$$
\dim \Gr^W_l \Hcts^2(\t_{g,r}^n)_\lambda \le
\dim Gr^W_{l+|\lambda|} H^2(\M_{g,r}^n,\V(\lambda)). \qed
$$
\end{corollary}

So, in order to bound the degrees of the relations in $\t_{g,r}^n$,
it suffices to give a bound the weights on $H^2(\M_{g,r}^n,\V(\lambda))$.
In the absolute case we have the following result of Kabanov
\cite{kabanov,kabanov:purity} which is proved using intersection homology.

\begin{theorem}[Kabanov]
For each irreducible rational representation $V(\lambda)$ of $Sp_g$,
we have
$$
\Gr^W_{k + |\lambda|}H^2(\M_g,\V(\lambda)) = 0
$$
when
$$
\begin{cases}
k \neq 2 & \text{ when $g\ge 6$;}\cr
k \neq 2,3 & \text{ when $3 \le g < 6$.}
\end{cases}
$$
\end{theorem}

Combining Kabanov's result with the previous results, we obtain:

\begin{corollary}
If $g\ge 3$, then $\Gr^W_\dot\t_g$ has a presentation with only
quadratic and cubic relations, and only quadratic relations when
$g\ge 6$. \qed
\end{corollary}

It is now an easy matter to insert the decorations:

\begin{corollary}
If $g\ge 3$, then $\Gr^W_\dot\t_{g,r}^n$ has a presentation with only
quadratic and cubic relations, and only quadratic relations when
$g\ge 6$. \qed
\end{corollary}

\section{The Lower Central Series Quotients of a Surface Group}

In this section we gather some information about $\Gr\p_g^1$
that will be useful when computing relations in $\Gr\t_g$ and
$\Gr\t_g^1$. Our basic tool, once again, is continuous cohomology.

A group is called {\it pseudo-nilpotent} if $\theta$ is an
isomorphism. A proof of the following result is sketched
by Kohno and Oda in \cite[(4.1)]{kohno-oda}.

\begin{theorem}\label{curve}
If $g \ge 1$, then $\pi_g^1$ is pseudo-nilpotent. \qed
\end{theorem}

Even though we will not be needing it, we record the following
result which is stated by Kohno and Oda \cite[(4.1)]{kohno-oda}.
Their proof is incorrect --- cf.\  (\ref{error}). Nonetheless, the
result follows directly from (\ref{curve}) and \cite[(5.7)]{hain:cycles}.

\begin{corollary}\label{kohno-oda}
If $g=0$ and $r\ge 1$, or if $g \ge 1$, then, for all $n\ge 0$,
each of the decorated pure braid groups $F_{g,r}^n$ is
pseudo-nilpotent. \qed
\end{corollary}

Since $H_1(\p_g^1)$ is the fundamental representation of $\sp_g$,
$\Gr^W_\dot \p_g^1$ is a graded Lie algebra in $R(\sp_g)$, and its
complex of chains $\Lambda^\dot \Gr^W_\dot\p_g^1$ is a complex in
$R(\sp_g)$.

We shall write $\p_g$ \label{p_def} instead of $\p_g^1$, and $\pi_g$
\label{pi_def} instead of $\pi_g^1$. We shall denote the $l$th weight
graded quotient of a Hodge Lie algebra $\g$ by \label{gr_def}
$\g(l)$. In particular, we shall denote $\Gr^W_{-l}\p_g^1$  by
\label{pgr_def} $\p_g(l)$.

\begin{corollary}\label{complex}
If $g\ge 1$, then, for each $l\ge 3$, the complex
$$
\Gr^W_{-l} \Lambda^\dot\Gr^W_\dot\p_g
$$
is an acyclic complex of $\sp_g$ modules. When $l=2$, we have an
exact sequence
$$
0 \to \Q(1) \to \Lambda^2 \p_g(1) \to \p_g(2) \to 0
$$
of $\sp_g$ modules. \qed
\end{corollary}

This result allows us to compute the $\p_g(l)$ inductively
as elements of $R(\sp_g)$. As before, we fix a set
$\lambda_1,\dots,\lambda_g$ of fundamental weights of $\sp_g$.

\begin{proposition}\label{lcs_quots}
For all $g\ge 3$, the highest weight decomposition of $\p_g(l)$
when $1\le l \le 4$ is given by
$$
\p_g(l) =
\begin{cases}
V(\lambda_1) & \text{ when $l=1$}; \cr
V(\lambda_2) & \text{ when $l=2$}; \cr
V(\lambda_1 + \lambda_2) & \text{ when $l=3$}; \cr
V(2\lambda_1) + V(2\lambda_1+\lambda_2) +
V(\lambda_1 + \lambda_3) & \text{ when $l=4$}.
\end{cases}
$$
\end{proposition}

\begin{proof}
This is a straightforward consequence of (\ref{complex}). To show
how this works, we prove the case where $l=3$. In this case, we
have the exact sequence
$$
0 \to \Lambda^3 \p_g(1) \to \p_g(1)\otimes \p_g(2) \to \p_g(3) \to 0
$$
in $R(\sp_g)$. Taking euler characteristics and applying the result
for $l=2$ and $g\ge 3$, we see that
$$
\p_g(3) =
V(\lambda_1)\otimes V(\lambda_2) - \Lambda^3 V(\lambda_1)
= V(\lambda_1 + \lambda_2).
$$
\end{proof}

Since the $k$th exterior power is the Schur functor corresponding
to the Young diagram with $k$ rows and one box in each row, and
since $\p_g(1) = H_1(\pi_g)$ is the fundamental representation
of $\sp_g$, we obtain the following stability result for the
graded quotients of the lower central series of $\pi_g$.

\begin{corollary}
The highest weight decomposition of $\p_g(l)$ is independent of
$g$ when $l\ge g$. \qed
\end{corollary}

\section{The Action of $\t_g^1$ on $\p_g$}
\label{inf_action}

In this section we obtain a lower bound for the size of
$\Gr^W_l\t_g$ when $l=2,3$ and $g\ge 3$ by studying the action of
$\t_g^1$ on $\p_g$. This will provide an upper bound on the size of
the quadratic and cubic relations of $\Gr^W_\dot\t_g$. Related
results have been obtained Asada-Kaneko \cite{japanese}, Morita
\cite{morita:trace} and Asada-Nakamura \cite{asada-nakamura}. We also
use a result \cite{asada-nakamura} of Asada and Nakamura to prove
that $\t_g$ is infinite dimensional. (This fact also follows from
a recent result of Oda \cite{oda}.)

First, some notation. Denote the pronilpotent Lie algebra
$W_{-1}\Der \p_g$ by \label{der_def} $\d_g$, and the quotient of this
by inner automorphisms by \label{out_def} $\o_g$. Once a base point $x$
of $\M_g^1$ has been chosen, each of these acquires the structure
of a Hodge Lie algebra.

We have natural representations
$$
\u_g^1 \to \d_g \text{ and } \u_g \to \o_g.
$$
These induce homomorphisms of their associated graded Lie algebras.

It is clear that there is an injective homomorphism
$$
\d_g(l) \hookrightarrow \Hom(\p_g(1),\p_g(l+1)).
$$
Each element $\delta : \p_g(1) \to \p_g(l+1)$ determines a derivation
$\tilde{\delta}$ of the free Lie algebra $\L(\p_g(1))$, the second
graded quotient of which is isomorphic to $\p_g(2) \oplus \Q\omega$,
where $\Q \omega$ is the unique copy of the trivial representation
in $\Lambda^2 \p_g(1)$. By taking the image of $\tilde{\delta}(\omega)$
under the projection
$$
\L(\p_g(1)) \to \Gr^W_\dot \p_g,
$$
we obtain an element $\sigma_g(\delta)$ of $\p_g(l+1)$. Observe
that $\delta$ induces a derivation of $\Gr\p_g$ if and only if
$\sigma_g(\delta)$ vanishes. We therefore have a surjection
$$
\d_g \to \ker\left\{\Hom(\p_g(1),\p_g(l+1))
\stackrel{\sigma}{\to} \p_g(l+2)\right\}.
$$

\begin{proposition}\label{gr_der}
The map $\sigma$ is surjective. Consequently,
$$
\d_g(l) = \p_g(1)\otimes\p_g(l+1) - \p_g(l+2)
$$
in $R(\sp_g)$.
\end{proposition}

\begin{proof}
Consider the diagram
$$
\begin{CD}
\Hom(\p_g(1),\p_g(l+1)) @>{\sigma_g}>> \p_g(l+2)\cr
@VVV @| \cr
\p_g(1)\otimes \p_g(l+1) @>{[\blank,\blank]}>> \p_g(l+2)
\end{CD}
$$
where the left hand vertical map is induced by the quadratic
form $\omega = \sum a_i\wedge b_i$. This diagram commutes as the
left hand map satisfies
$$
\delta \mapsto
\sum_{i=1^g} a_i\otimes\delta(b_i) - b_i\otimes\delta(a_i),
$$
which goes to
$$
\sigma_g(q) =
\sum_{i=1}^g\left([\delta(a_i),b_i] + [a_i,\delta(b_i)]\right)
$$
under the bracket. Since the bottom map is surjective and all
maps are $\sp_g$ equivariant, the result follows.
\end{proof}

Combining this with the computation of the first few graded
quotients of $\p_g$ given in (\ref{lcs_quots}), we obtain the following
result.

\begin{corollary}\label{der_quots}
For all $g\ge 3$, we have
$$
\d_g(l) =
\begin{cases}
V(\lambda_3)+V(\lambda_1) &\text{ when }l=1;\cr
V(2\lambda_2)+V(\lambda_2)&\text{ when }l=2;\cr
V(2\lambda_1+\lambda_3)+V(\lambda_1+\lambda_2)+V(3\lambda_1)&
\text{ when } l=3.
\end{cases}
$$ \qed
\end{corollary}

The computation of $\d_g(1)$ is simply another formulation of the
Johnson homomorphism.

It is proven in \cite[$A^\prime$, p.~149]{japanese} that the center
of $\Gr\p_g$ is  trivial, so that the inclusion $\p_g \to \d_g$ of
the inner automorphisms is injective.

\begin{proposition}
For all $g\ge 3$ and all $l\ge 1$,
$\o_g(l) = \d_g(l) - \p_g(l)$. \qed
\end{proposition}

Combining (\ref{lcs_quots}) and (\ref{der_quots}), we obtain
the following computation.

\begin{corollary}\label{out_quots}
For all $g\ge 3$, we have
$$
\o_g(l) =
\begin{cases}
V(\lambda_3) &\text{ when }l=1;\cr
V(2\lambda_2)&\text{ when }l=2;\cr
V(2\lambda_1+\lambda_3)+V(3\lambda_1)&
\text{ when } l=3.
\end{cases}
$$ \qed
\end{corollary}

It does not seem obvious {\it a priori}, that $\t_g$ is infinite
dimensional.\footnote{This result also follows quite directly from
a result of Oda \cite{oda}.}

\begin{proposition}\label{quotients}
For all $g\ge 3$, the image of $\t_g$ in $\o_g$ is infinite dimensional.
\end{proposition}

\begin{proof}
Since $\t_g \to \o_g$ is a morphism of MHS, the image $\g$ has a MHS.
Since $\Gr^W_\dot$ is an exact functor,
$$
\Gr^W_\dot\g = \text{ image of }\{\Gr^W_\dot\t_g \to \Gr^W_\dot\o_g\}.
$$
So it suffices to show that each graded quotient of $\g$ is non-trivial.
It follows from the result Asada and Nakamura
\cite[Theorem~B]{asada-nakamura} that the image of
$$
\t_g^1(l) \to \d_g(l)
$$
contains the representation $V(2m\lambda_1 + \lambda_3)$ when
$l = 2m+1$, and $V(2m\lambda_1 + 2\lambda_2)$ when $l = 2m+2$.
These representations both have the maximal possible depth, $l+2$.
But the inner automorphisms $\p_g(l)$ in $\d_g(l)$
have depth at most $l$. The result follows.
\end{proof}

We can now bound below the low degree relations in $\t_g$.

\begin{proposition}\label{upper}
For all $g\ge 3$, the image of $\u_g(l)$ in $\o_g(l)$ is
$$
\begin{cases}
V(\lambda_3) & \text{ when } l=1;\cr
V(2\lambda_2) & \text{ when } l=2;\cr
V(2\lambda_1+\lambda_3) & \text{ when } l=3.
\end{cases}
$$
\end{proposition}

\begin{proof}
It follows from (\ref{quotients}) that when $g\ge 3$, the image of
$\u_g(l) \to \o_g(l)$ is non-trivial for all $l$. Since this map is
$\sp_g$ equivariant, the image of $\u_g(2)$ must be all of $\o_g(2)$.
Since $\Gr^W_\dot \u_g$ is  generated by $\u_g(1)$, and since
$V(3\lambda_1)$ does not appear in $\u_g(1)\otimes \u_g(2)$, the
assertion for $l=3$ follows.
\end{proof}

Note that the copy of $V(3\lambda_1)$ is detected by Morita's
trace \cite{morita:trace}.

\section{Quadratic Relations}
\label{quadratic_relns}

In this section, we find some obvious quadratic relations in $\t_g$
for each $g\ge 4$. These give a lower bound for the relations in $\t_g$.
Serendipitously, this coincides with the upper bound (\ref{upper})
derived in the previous section, thus yielding all the quadratic
relations.

\begin{theorem}\label{lower}
For all $g \ge 3$, we have
$$
\Gr^W_{-2} \t_g = \Gr^W_{-2} \u_g = V(2\lambda_2) + V(0).
$$
\end{theorem}

The proof occupies the rest of this section. We prove the result by
finding a pair of commuting elements  $\phi$ and $\psi$ of the $T_g$
whose class
$$
\tau(\phi)\wedge\tau(\phi) \in \Lambda^2 V(\lambda_3)
$$
generates the $Sp_g$ complement of $V(2\lambda_2) + V(0)$ for all
$g \ge 3$. Since we know, by (\ref{upper}), that the quadratic relations
are contained in the complement of $V(2\lambda_2) + V(0)$, we have found
all quadratic relations.

\begin{lemma}\label{computation}
If $g\ge 3$, then
\begin{multline*}
\Lambda^2 V(\lambda_3) = \\
\begin{cases}
V(\lambda_6) + V(\lambda_4) + V(\lambda_2) + V(\lambda_2 + \lambda_4)
+ V(2\lambda_2) + V(0) & \text{when $g\ge 6$;}\cr
V(\lambda_4) + V(\lambda_2) + V(\lambda_2 + \lambda_4)
+ V(2\lambda_2) + V(0) & \text{when $g=5$;}\cr
V(\lambda_2) + V(\lambda_2 + \lambda_4) + V(2\lambda_2) + V(0) &
\text{when $g = 4$;}\cr
V(2\lambda_2) + V(0) & \text{when $g=3$.}
\end{cases}
\end{multline*}
\qed
\end{lemma}

{}From (\ref{mhs_torelli}), we know that $V(0)$ occurs in $\t_g(2)$.
By (\ref{upper}) and the previous proposition, there is nothing to
prove when $g=3$. So we suppose that $g\ge 4$.

We use the notation introduced in Section \ref{reps}. Set
$$
\omega = a_1\wedge b_1 + \cdots + a_g\wedge b_g.
$$

\begin{proposition}\label{construction}
When $g\ge 3$, there are elements $\phi_{i,j}$, $1\le i < j \le g$
of the Torelli group whose image under the Johnson homomorphism
$$
\tau_g : H_1(T_g) \to V(\lambda_3)
$$
is given by
$$
(g-1)\tau_g(\phi_{i,j}) = (g-1) a_i\wedge a_j
\wedge b_j - a_i \wedge \omega
$$
Here we are viewing $V(\lambda_3)$ as a submodule of
$\Lambda^3 V(\lambda_1)$.  Moreover, we can choose them such that
$\phi_{1,2}$ and $\phi_{3,4}$ commute when $g\ge 4$.
\end{proposition}

\begin{proof}
For $1 \le i < j \le g$ is easy to construct elements $\phi_{i,j}$
of the pointed Torelli group $T^1_g$ with
$$
\tau_g^1(\phi_{i,j}) = a_i\wedge a_j \wedge b_j \in \Lambda^3
V(\lambda_1).
$$ (To compute $\tau_g^1: H_1(T_g^1)\to \Lambda^3V(\lambda_1)$, use
Johnson's original
definition in terms of the action of $T_g^1$ on $\pi_g$.)
\vspace*{1.5in}\\
It is also easy to arrange for $\phi_{1,2}$ and $\phi_{3,4}$ to have
disjoint supports, and therefore commute.
To compute $\tau_g(\phi_{i,j}) \in V(\lambda_3)$, we just use the fact
that the maps
$$
\underline{\blank}\wedge\omega : V(\lambda_1)
\to \Lambda^3 V(\lambda_1)
$$
and
$$
p: \Lambda^3 V(\lambda_1) \to V(\lambda_1)
$$
defined by $p(x\wedge y\wedge z) = q(x,y)z + q(y,z) x + q(z,x) y$ are
$\sp_g$-equivariant and satisfy
$$
p\circ (\underline{\blank}\wedge\omega) = (g-1)\id.
$$
It follows that $V(\lambda_3)$ is the kernel of $p$ and that
$$
(g-1) \tau_g(\phi_{i,j}) = (g-1) a_i\wedge a_j\wedge b_j - a_i \wedge
\omega.
$$
\end{proof}

Take $\phi = \phi_{1,2}$ and $\psi = \phi_{3,4}$. Since these commute,
$$
v :=\tau_g(\phi)\wedge\tau_g(\psi) \in \Lambda^2 V(\lambda_3)
$$
will lie inside the $\sp_g$ module of quadratic relations. Denote
the $\sp_g$ submodule of $\Lambda^2 V(\lambda_3)$ generated by $v$
by $V$. By (\ref{construction}),
$$
v =
[(g-1)\, a_1\wedge a_2 \wedge b_2 - a_1\wedge \omega]
\wedge
[(g-1)\, a_3\wedge a_4 \wedge b_4 - a_3\wedge \omega].
$$
Recall that elements of $\sp_g$ act on exterior powers as derivations.
Note also that for all $X \in \sp_g$, $X\omega = 0$. We now compute
the highest weight decomposition of $V$.
\smallskip

\noindent{$\mathbf \lambda_2 + \lambda_4$:} Apply $F_{2,3}$, then
$F_{1,4}$, then $T_{2,3}$ to $v$ to get
$$
(g-1)^2
[a_1\wedge a_2 \wedge a_3 ] \wedge [a_1\wedge a_2 \wedge a_4]
\in \Lambda^2\Lambda^3 V(\lambda_1)
$$
which is a highest weight vector on which $\h$ acts via the character
$$
\lambda_2 + \lambda_4 = (t_1 + t_2) + ( t_1 + t_2 + t_3 + t_4).
$$

To decompose the rest of $V$, consider the $\sp_g$-equivariant map
$$
\Lambda^2 V(\lambda_3) \hookrightarrow \Lambda^2\Lambda^3
V(\lambda_1) \to \Lambda^6 V(\lambda_1).
$$
Denote the image of $V$ under this map by $W$. It is spanned by
the image of $v$ in $\Lambda^6 V(\lambda_1)$. For the time being,
we suppose that $g\ge 6$.
\smallskip

\noindent{$\mathbf \lambda_6$:} In this case, the image of $v$ in $W$
is
$$
w := ((g-1)\, a_1\wedge a_2 \wedge b_2 - a_1\wedge \omega)
\wedge
((g-1)\, a_3\wedge a_4 \wedge b_4 - a_3\wedge \omega).
$$
To find a highest weight vector for the representation it generates,
first apply $F_{2,5}$, then $F_{4,6}$ to this vector to get the highest
weight vector
$$
(g-1)^2
a_1\wedge a_2 \wedge a_3 \wedge a_4 \wedge a_5 \wedge a_6
$$
of $W$ on which $\h$ acts via the character
$$
\lambda_6 = t_1 + t_2 + t_3 + t_4 + t_5 + t_ 6.
$$

To show that the weights $\lambda_2$ and $\lambda_4$ occur in
$V$, it is useful to recall that for all $k\ge 2$, there is an $\sp_g$
equivariant projection
\begin{equation}\label{contraction}
\theta_k : \Lambda^k V(\lambda_1) \to \Lambda^{k-2}
V(\lambda_1)
\end{equation}
which is defined by
$$
x_1 \wedge \ldots \wedge x_k \mapsto \sum_{1 \le i < j \le k}
(-1)^{i+j+1} q(x_i,x_j)\, x_1 \wedge \ldots \wedge
\hat{x_i}\wedge  \ldots \wedge \hat{x_j} \wedge \ldots
\wedge x_k.
$$

\noindent{$\mathbf \lambda_4$: } The image of $V$  in $\Lambda^4
V(\lambda_1)$ is generated by $\theta_6(w)$ which is
$$
(g-1)^2
a_1 \wedge a_3 \wedge (a_2 \wedge b_2 + a_4 \wedge b_4)
- (g-1)(g-3)\,
a_1 \wedge a_3 \wedge(a_2 \wedge b_2 + a_4 \wedge b_4)
$$ $$
- 2(g-1)\, a_1 \wedge a_3 \wedge \omega
 - 2(g-2)\, a_1 \wedge a_3 \wedge \omega
$$ $$
= 2(g-1)\, a_1 \wedge a_2 \wedge (a_3 \wedge b_3 + a_4 \wedge
b_4) -2\, a_1 \wedge a_2 \wedge \omega.
$$
Applying $F_{3,6}$, then $S_{4,6}$, one gets the highest weight
vector
$$
2(g-1)\, a_1\wedge a_2\wedge a_3\wedge a_4
$$
on which $\h$ acts via the character
$\lambda_4 = t_1 + t_2 + t_3 + t_4$.
\smallskip

\noindent{ $\mathbf \lambda_2$:} The image of $V$ in $\Lambda^2
V(\lambda_1)$ is generated by the image under $\theta_4$ of
$\theta_6(w)$.  This is
$$
4(g-1)\, a_1 \wedge a_3 - 2(g-2)\, a_1 \wedge a_3 = 2g\, a_1 \wedge
a_3.
$$
Apply $S_{2,3}$ to this to get $2g\, a_1 \wedge a_2$  upon which
$\h$ acts with highest weight $\lambda_2 = t_1 + t_2$.

We sketch the remaining cases $g=4,5$. When $g=5$, $W$ is generated
by the vector
$$
(g-1)\left(a_1\wedge a_2\wedge a_3\wedge b_2 -
a_1\wedge a_3\wedge a_4\wedge b_4\right)\wedge \omega.
$$
By contracting with $q$ as above, it is easy to see that this
vector generates a submodule of
$$
\omega \wedge \Lambda^4 V(\lambda_1) \cong \Lambda^6 V(\lambda_1)
$$
isomorphic to $V(\lambda_4) + V(\lambda_2)$.

When $g=4$, $W$ is generated by $a_1\wedge a_3 \wedge \omega^2$.
Again, by contracting with $q$, it is easy to see that this vector
generates a submodule of
$$
w^2\wedge \Lambda^2 V(\lambda_1)
\subset \Lambda^6 V(\lambda_1)
$$
isomorphic to $V(\lambda_2)$.

\begin{remark}
Note that we have determined the quadratic relations for all
$g\ge 3$. One should be able to determine the cubic relations when
$3 \le g \le 5$ by applying similar methods and the fact that
the Dehn twist about the separating curve $C$ below commutes
with the bounding pair map associated to the curves $C'$ and
$C''$. The Dehn twist about $C$ is in the kernel of the Johnson
homomorphism, but has non-trivial image in the second graded
quotient of $\t_g$.
\vspace*{2in}
Note that there have to be cubic relations in genus 3
as there are no quadratic relations, and there has to be one copy
of $V(\lambda_3)$ in the cubic relations to ensure the existence of
the central $\Ga$.
\end{remark}

\section{Presentations of $\t_g$, $t_{g,1}$ and $\t_g^1$}
\label{special}

Recall that $\L(V)$ denotes the free Lie algebra generated by the
vector space $V$. In this section we shall give presentations of
$\t_g$, $\t_g^1$ and $\t_{g,1}$ when $g\ge 6$. First $\t_g$ ---
combining (\ref{upper}), (\ref{lower}) and (\ref{construction}), we
have:

\begin{theorem}
For all $g\ge 6$, $\Gr^W_\dot \t_g$ is isomorphic to
$$
\L(V(\lambda_3))/R_g
$$
as a graded Lie algebra in $R(\sp_g)$, where $R_g$ is the ideal
generated by the quadratic relations
$$
V(\lambda_6) + V(\lambda_4) + V(\lambda_2) + V(\lambda_2 +\lambda_4)
\subseteq \Lambda^2 V(\lambda_3). \qed
$$
\end{theorem}

Since $\t_g\otimes\C \cong \prod_l \Gr^W_l \t_g\otimes\C$, this gives
the desired presentation of $\t_g$ for $g\ge 6$.

We next consider $\t_g^1$.
Fix a point $[C;x]$ of $\M_g^1$ so that $\t_g^1$ and $\p_g$
have canonical MHSs. Then the sequence
$$
0 \to \p_g^1 \to \t_g^1 \to \t_g \to 0
$$
is an exact sequence of MHSs. Since $\Gr^W_\dot$ is an exact functor,
the sequence
$$
0 \to \Gr^W_\dot \p_g^1 \to \Gr^W_\dot \t_g^1 \to \Gr^W_\dot \t_g \to 0
$$
is exact in $R(\sp_g)$.  Since the sequence of $H_1$'s is
canonically split in $R(\sp_g)$, there is a canonical lift
$$
\L(H_1(\t_g)) \to \Gr^W_\dot\t_g^1
$$
of the natural homomorphism $\L(H_1(\t_g)) \to \Gr^W_\dot\t_g$. Since
$\p_g(2)$ is isomorphic to $V(\lambda_2)$, and since this representation
does not occur in $\t_g(2)$, we can take the $\lambda_2$ component of the
bracket
$$
\Lambda^2 H_1(\t_g) \hookrightarrow \Lambda^2 H_1(t_g^1) \to \t_g^1(2)
$$
to obtain an $\sp_g$ module map
\begin{equation}\label{bracket}
c : \Lambda^2 H_1(\t_g) \to \p_g(2) \cong V(\lambda_2).
\end{equation}
This and the map
\begin{equation}\label{action2}
H_1(\t_g) \otimes H_1(\p_g) \to \p_g(2)
\end{equation}
induced by the bracket completely determine $\Gr^W_\dot \t_g^1$ given
$\Gr^W_\dot \t_g$ and $\Gr^W_\dot \p_g$. The map (\ref{action2})
is simply the adjoint of the Johnson homomorphism
$$
\tau_g^1 : H_1(\t_g) \to \Lambda^3 V \subset \Hom(H_1(\p_g),\p_g(2)).
$$
So, to give a presentation of $\t_g^1$, we have to determine the map
(\ref{bracket}). We do this by studying the action of $\L(V(\lambda_3))$
on $\L(V(\lambda_1))$.

Set $V=V(\lambda_1)$. We identify $V$ with $H_1(C)$ and $H_1(T_g^1)$
with $\Lambda^3 V$ via the Johnson homomorphism. Recall that $V(\lambda_3)$
can be realized as the kernel of the map $p:\Lambda^3 V \to V$ defined by
\begin{equation}\label{projn}
p : v_1\wedge v_2 \wedge v_3 \mapsto (v_1\cdot v_2)v_3 +
(v_2\cdot v_3) v_1 + (v_3\cdot v_1) v_2.
\end{equation}
We identify $H_1(T_g) \cong \Lambda^3 V/V$ with $V(\lambda_3)$
via the map
$$
V(\lambda_3) = \ker p \hookrightarrow \Lambda^3 V \to \Lambda^3 V/V.
$$

The natural action of $\Lambda^3 V$ on $\L(V)$ is defined by
\begin{multline*}
e_1\wedge e_2 \wedge e_3 \mapsto -\left\{v \mapsto (e_1\cdot v)
[e_2,e_3] + (e_2\cdot v)[e_3,e_1] + (e_3\cdot v)[e_1,e_2]\right\} \cr
\in \Hom(V,\Lambda^2 V) \subseteq \Der \L(V).
\end{multline*}
(With this choice of sign, $\sum x\wedge a_j\wedge b_j\mapsto\ad(x)$.)
It follows from the definition of the Johnson homomorphism that the
composite
$$
H_1(T_{g,1}) \stackrel{\tau_{g,1}}{\longrightarrow}
\Lambda^3 V \hookrightarrow \Hom(V,\Lambda^2 V)
$$
is the map induced by the action of $T_{g,1}$ on $\pi_{g,1}$.
The action descends to the action of $\Gr^W_\dot\t_g^1$ on
$\Gr^W_\dot \p_g$.

Define a projection $r : \Lambda^2 V(\lambda_3) \to V(\lambda_2)$
to be the composite
$$
\Lambda^2 V(\lambda_3) \hookrightarrow \Lambda^2 \Lambda^3 V
\stackrel{\text{mult}}{\longrightarrow} \Lambda^6 V
\stackrel{\theta_4\theta_6}{\longrightarrow} \Lambda^2 V
\to V(\lambda_2)
$$
where $\theta_k$ is the contraction (\ref{contraction}) defined in
the previous section, and the last map is the standard projection
$$
u\wedge v \mapsto u\wedge v - (u\cdot v)\,\omega/(g-1).
$$
Since there is only one copy of $V(\lambda_2)$
in $\Lambda^2 V(\lambda_3)$, this projection is unique up to a scalar.

\begin{proposition}\label{bra_const}
The map (\ref{bracket}) is given by
$$
c[u,v] = -\frac{1}{2g+2} \ad(r(u\wedge v)) \in \Hom(H_1(p_g),\p_g(3)).
$$
In particular, this map is non-zero, and the extensions
$$
0 \to \p_g \to \t_g^1 \to \t_g \to 0 \text{ and }
0 \to \p_g \to \u_g^1 \to \u_g \to 0
$$
are not split.
\end{proposition}

We now sketch the proof.
Denote the degree $k$ part of the free Lie algebra $\L(V)$ by
$\L(V)(k)$.
Recall that there is a standard exact sequence
$$
0 \to \Lambda^3 V \stackrel{j}{\to} V\otimes \Lambda^2 V
\stackrel{b}{\to} \L(V)(3) \to 0.
$$
The first map is the ``Jacobi identity'' map
$$
j : v_1\wedge v_2 \wedge v_3 \mapsto v_1\otimes v_2\wedge v_3
+ v_2 \otimes v_3\wedge v_1 + v_3 \otimes v_1 \wedge v_2,
$$
and the second map is the bracket. (We are identifying $\L(V)(2)$
with $\Lambda^2 V$ in the standard way.)

\begin{lemma}\label{compn}
The bracket $[e_1\wedge e_2 \wedge e_3, f_1\wedge f_2 \wedge f_3]$
of two elements of $\Lambda^3 V$ as derivations of $\L(V)$ is obtained
by summing the expression
\begin{multline*}
(e_1\cdot f_1)
\big(
e_2\otimes [e_3,[f_2,f_3]] - e_3\otimes [e_2,[f_2,f_3]]
+ f_2\otimes [f_3,[e_2,e_3]] - f_3\otimes [f_2,[e_2,e_3]]
\big) \cr
\quad \in V\otimes \L(V)(3) \cong \Hom(V,\L(V)(3))
\end{multline*}
cyclically in $(e_1,e_2,e_3)$ and in $(f_1,f_2,f_3)$. \qed
\end{lemma}

We shall view this expression as an element of
$\left(V\otimes V \otimes \Lambda^2 V\right)
/\left(V \otimes \Lambda^3 V\right)$.
The next step is to write down the projections of this group onto
$V(\lambda_2)$.

There are four copies of $V(\lambda_2)$ in
$V\otimes V \otimes \Lambda^2 V$. These are detected by the
following four projections onto $\Lambda^2 V$:
\begin{align*}
p_1 : u_1\otimes u_2 \otimes u_3 \wedge u_4
&\mapsto (u_1\cdot u_2)u_3\wedge u_4 \cr
p_2 : u_1\otimes u_2 \otimes u_3 \wedge u_4
&\mapsto (u_3\cdot u_4) u_1\wedge u_2 \cr
p_3 : u_1\otimes u_2 \otimes u_3 \wedge u_4
&\mapsto
\big((u_1\cdot u_4)u_2\wedge u_3 - (u_1\cdot u_3)u_2\wedge u_4\big)/2\cr
p_4 : u_1\otimes u_2 \otimes u_3 \wedge u_4
&\mapsto
\big((u_2\cdot u_3)u_1\wedge u_4 - (u_2\cdot u_4)u_1\wedge u_3\big)/2.
\end{align*}
One can easily check that there are two copies of $V(\lambda_2)$ in
$V\otimes \Lambda^3 V$ and that the projections $p_1 - p_3$ and
$p_2 - p_4$ vanish on these. This leaves two copies of $V(\lambda_2)$
in $\left(V\otimes V \otimes \Lambda^2 V\right)
/\left(V \otimes \Lambda^3 V\right)$. One of these vanishes in
$\Hom(V,\p_g(3))$ as $V\otimes V \otimes \omega$ projects to zero there.
We are now ready to compute.

Since
$$
u_j = a_j\wedge a_3 \wedge b_3 - a_j\wedge a_4 \wedge b_4
$$
lies in the kernel of the projection $p$ above when $j=1,2$,
$u_1\wedge u_2$ is an element of $\Lambda^2 V(\lambda_3)$. The
projection $r$ takes $u_1\wedge u_2$ to $-4\, a_1\wedge a_2$.
On the other hand, by straightforward computations using
(\ref{compn}), we have
$$
p_1([u_1,u_2]) = p_2([u_1,u_2]) = 0,\text{ and }
p_3([u_1,u_2]) = p_4([u_1,u_2]) = -4\, a_1\wedge a_2.
$$
Consequently,
$$
(p_1 - p_3)([u_1,u_2]) = (p_2 - p_4)([u_1,u_2]) = 4\, a_1\wedge a_2.
$$

Next observe that $\ad [a_1,a_2]$ corresponds to the element
$$
- \sum_{j=1}^g \left( a_j\otimes b_j \otimes a_1\wedge a_2
- b_j \otimes a_j \otimes a_1 \wedge a_2\right).
$$
of $\left(V\otimes V \otimes \Lambda^2 V\right)
/\left(V \otimes \Lambda^3 V\right)$.
Since $\sum\, [a_j,b_j] = 0$, $\ad [a_1,a_2]$ is also represented
by
$$
z := - \sum_{j=1}^g \left( a_j\otimes b_j \otimes a_1\wedge a_2
- b_j \otimes a_j \otimes a_1 \wedge a_2\right)
- 2 \sum_{j=1}^g a_1\otimes a_2\otimes a_j\wedge b_j.
$$
By direct computation, we have
$$
(p_1-p_3)(z) = (p_2-p_4)(z) = -(2g+2)\, a_1\wedge a_2.
$$\,
This concludes the proof of Proposition \ref{bra_const}. \qed

Next we consider the case of $\t_{g,1}$. Fix a point $(C;x,v)$ of
$\M_{g,1}$ so that $\t_{g,1}$, $\p_{g,1}$, etc.\ all have compatible
MHSs. By strictness, the sequence
$$
0 \to \Gr^W_\dot \p_{g,1} \to \Gr^W_\dot \t_{g,1}
\to \Gr^W_\dot \t_g \to 0
$$
is exact in $R(\sp_g)$. Since the sequence
$$
0 \to \Q(1) \to \p_{g,1} \to \p_g \to 0
$$
is exact, it follows that $\p_{g,1}(2)$ is isomorphic to $\Lambda^2 V$
via the bracket.

As in the case of $\t_g^1$, there is a canonical lifting
$\L(H_1(\t_g)) \to \Gr^W_\dot\t_{g,1}$
of the natural surjection $\L(H_1(\t_g)) \to \Gr^W_\dot \t_g$.
It follows that to give a presentation of $\Gr^W_\dot \t_{g,1}$ given
presentations of $\t_g$ and $\p_{g,1}$, it suffices to give the map
\begin{equation}\label{first}
H_1(\t_g)\otimes H_1(\p_{g,1}) \longrightarrow \p_{g,1}(2)
\end{equation}
induced by the bracket, together with the $\lambda_2$ component
\begin{equation}\label{second}
\Lambda^2 H_1(\t_g) \longrightarrow V(\lambda_2) \subset \p_{g,2}
\end{equation}
and the invariant part
\begin{equation}\label{third}
c_0 : \Lambda^2 H_1(\t_g) \longrightarrow
\t_{g,1}(2)^{Sp_g} \cong \Q(1)^2
\end{equation}
of the bracket. As in the case of $\t_g^1$, the first map (\ref{first})
is the adjoint of the Johnson homomorphism and the second (\ref{second}),
by naturality with respect to the projection $\t_{g,1} \to \t_g^1$, is
the map $c$ determined in (\ref{bra_const}). It remains to determine the
map (\ref{third}). Observe that the sequence
$$
0 \to \p_{g,1}(2)^{Sp_g} \to \t_{g,1}(2)^{Sp_g} \to \t_g(2)^{Sp_g} \to 0
$$
splits canonically as the canonical central $\Ga$ in $\t_{g,1}$
projects to the canonical central $\Ga$ in $\t_g$ by (\ref{central_ext}),
and because $\Ga = \t_g(2)^{Sp_g}$. As a generator of $\p_{g,1}(2)^{Sp_g}$
we take $\sum\, [a_j,b_j]$.

Fix an invariant bilinear form $\bil \blank \blank$ on $V(\lambda_3)$ by
insisting that
$$
\bil {a_1\wedge a_2 \wedge a_3} {b_1 \wedge b_2 \wedge b_3} = 1.
$$
We can therefore choose a generator $\gamma$ of $\Ga$ such that
if $u,v\in H_1(\t_g)$, then the invariant component of $[u,v]$
in $\t_g(2)$ is $\bil u v \, \gamma$.

\begin{proposition}\label{triv_cpt}
If $u,v \in H_1(\t_g)$, then
$$
c_0[u,v] =
\bil u v \,\gamma -
\frac{6\bil u v}{g(2g+1)}\sum_{j=1}^g\, [a_j,b_j].
$$
\end{proposition}

As in the previous case, we determine the coefficient by studying
the action of $\L(V(\lambda_3))$ on $\L(V)$. Note that $\Gamma_{g,1}$
acts on the free group $\pi_1(C - \{x\},v)$.%
\footnote{This notation denotes Deligne's fundamental group
of $C - \{x\}$ with base point the tangent vector $v \in T_x C$.} We
therefore have a representation $\t_{g,1} \to \p(C-\{x\},v) \cong
\L(V)$.\footnote{If we put the limit
MHS on $\pi_1(C - \{x\},v)$ associated with the tangent vector $v$,
then this action is a morphism of MHS.} We continue with the notation
in the proof of (\ref{bra_const}).

There are two copies of the trivial representation in
$V\otimes V \otimes \Lambda^2 V$. The corresponding projections
to $\Q$ are:
\begin{align*}
q_1 : u_1\otimes u_2 \otimes u_3 \wedge u_4
&\mapsto (u_1\cdot u_2)(u_3\cdot u_4) \cr
q_2 : u_1\otimes u_2 \otimes u_3 \wedge u_4
&\mapsto \bigl((u_1\cdot u_4)(u_2\cdot u_3)
- (u_1\cdot u_3) (u_2\cdot u_4)\bigr)/2.
\end{align*}
There is one copy of the trivial representation in $V\otimes
\Lambda^3 V$ and $q_1 - q_2$ vanishes on it. The vectors
$$
u_1 = a_1\wedge a_2 \wedge a_3
\text{ and } u_2 = b_1\wedge b_2 \wedge b_3
$$
both lie in $V(\lambda_3)$ and $\bil {u_1} {u_2} = 1$. It follows
from the formula (\ref{bracket}) that $[u_1,u_2]$ is obtained by
summing the expression
$$
a_2 \otimes [a_3,[b_2,b_3]] - a_3\otimes [a_2,[b_2,b_3]]
+ b_2 \otimes [b_3,[a_2,a_3]] - b_3 \otimes [b_2,[a_2,a_3]]
$$
over the cyclic group generated by $(1,2,3)$. We have
$(q_1 - q_2) ([u_1,u_2]) = 6$.
On the other hand, $\ad \sum[a_j,b_j]$ is represented by
$$
\sum_{j=1}^g \sum_{k=1}^g
(b_j\otimes a_j - a_j \otimes b_j)\otimes a_k\wedge b_k
$$
The projection $q_1 - q_2$ takes the value $-g(2g +1)$ on
this. The result follows.

\begin{remark}
The formulas in (\ref{bra_const}) and (\ref{triv_cpt}) are closely
related to those in Theorem~3.1 of Morita's paper
\cite{morita:cocycles}.
\end{remark}

\section{A Presentation of $\p_{g,r}^n$}
\label{braids2}

In this section we give an explicit quadratic presentation of the
pure braid Lie algebras $\p_{g,r}^n$ for all $g > 0$. We continue with
the notation of Section~\ref{braids1}. We fix a complex structure on
and a base point of $F_{g,r}^n$ by choosing a point
$$
[C;x_1,\dots,x_n;v_1,\dots,v_r]
$$
of $\M_{g,r}^n$.

In Section~\ref{braids1} we showed that $H^1(F_{g,r}^n(C))$ is pure of
weight 1. We will show that $H^2(F_{g,r}^n(C))$ is pure of weight 2,
from which the existence of a quadratic presentation of $\p_{g,r}^n$
will follow via Morgan's Theorem (\ref{morgan}).

First, some notation. Denote the projection
of $C^n$ onto its $i$th factor by $p_i$. Denote the image of the
inclusion
$$
p_i^\ast : H^\dot(C) \hookrightarrow H^\dot(C^n)
$$
by $H^\dot(C_i)$. For $x\in H^\dot(C)$, denote $p_i^\ast x$ by $x\sup{i}$.
Denote the component of $\Delta$ where the $i$th and $j$th coordinates
are equal by $\Delta_{ij}$. Fix a
symplectic basis $a_1, \dots, a_g,b_1,\dots, b_g$ of $H_1(C)$,
and let $\alpha_1,\dots, \alpha_g,\beta_1,\dots,\beta_g$ be the
dual basis of $H^1(C)$. Denote the positive integral generator of
$H^2(C)$ by $\zeta$, and the intersection form
$$
\sum_{r=1}^g \alpha_r \wedge \beta_r
$$
by $q$. When $i\neq j$, set
$$
q_{ij} = \sum_{r=1}^g \alpha_r\sup{i} \wedge \beta_r\sup{j}
+ \alpha_r\sup{j} \wedge \beta_r\sup{i}.
$$

\begin{lemma}
The Poincar\'e dual $PD(\Delta_{ij})$ of $\Delta_{ij}$ is
$\zeta\sup{i} + \zeta\sup{j} - q_{ij}$. \qed
\end{lemma}

This is elementary. Another elementary fact we shall need is the
following statement. It is easily proved using a Mayer-Vietoris
argument.

\begin{lemma}\label{isom}
The natural map
$$
\bigoplus_{i<j} H_{2n-3}(\Delta_{ij}) \to H_{2n-3}(\Delta)
$$
is an isomorphism. \qed
\end{lemma}

We can therefore write the Gysin map
$\gamma : H_{2n-3}(\Delta) \to H^3(C^n)$
as the sum of the Gysin maps
$\gamma_{ij} : H^1(\Delta_{ij}) \to H^3(C^n)$;
the map $\gamma_{ij}$ being given by cup product with $PD(\Delta_{ij})$.

\begin{lemma}\label{formula}
The composite
$H^1(C) \stackrel{p_k^\ast}{\longrightarrow} H^1(\Delta_{ij})
\stackrel{\gamma_{ij}}{\longrightarrow} H^3(C^n)$
is given by
$$
x \mapsto
\begin{cases}
\zeta\sup{i}\wedge x\sup{j} + \zeta\sup{j}\wedge x\sup{i} &
\text{ if $k \in \{i,j\}$;}\cr
\zeta\sup{i}\wedge x\sup{k} + \zeta\sup{j}\wedge x\sup{k}
- q_{ij}\wedge x\sup{k} & \text{ if $k\not\in \{i,j\}$. \qed}
\end{cases}
$$
\end{lemma}

It follows from
(\ref{h1_braid}) that the part
$$
0 \to \bigoplus_{i<j} \Z \to H^2(C^n) \to H^2(C^n - \Delta)
\to H_{2n-3}(\Delta) \to H^3(C^n)
$$
of the Gysin sequence is exact. Purity of $H^2(F_g^n(C))$ therefore
follows from the following proposition.

\begin{proposition}
The Gysin map $\gamma : H_{2n-3}(\Delta) \to H^3(C^n)$ is injective.
\end{proposition}

\begin{proof}
The Gysin sequence can be viewed as the fiber over $[C] \in \M_g$
of an exact sequence of (orbifold) local systems over $\M_g$. It
follows from (\ref{isom}) that the the monodromy actions of the last
two terms of the Gysin sequence above factor through the symplectic
group. It is  convenient, though not necessary, to decompose these
groups under the action of $Sp_g$.

First note that $H_{2n-3}(\Delta_{ij})$ is isomorphic to
$H^1(\Delta_{ij})$, which is isomorphic to $n-1$ copies of the
fundamental representation $V$. Next,
$$
H^3(C^n) = \bigoplus_{i\neq j} \left(H^2(C_i)\otimes H^1(C_j)\right)
\oplus \bigoplus_{i<j<k} H^1(C_i)\otimes H^1(C_j)\otimes H^1(C_k).
$$
Each of the terms $H^2(C_i)\otimes H^1(C_j)$ is a copy of the
fundamental representation that we shall denote by $V^i_j$. The
term $H^1(C_i)\otimes H^1(C_j)\otimes H^1(C_k)$ is isomorphic to
$V^{\otimes 3}$. It contains 3 copies of $V$. If, for $i,j,k$ distinct,
we set
$$
V_{ij}^k =
\text{ the image of }
\left\{H^1(C_k) \stackrel{\wedge q_{ij}}{\to} H^3(C^n)\right\},
$$
then
$$
\left[H^1(C_i)\otimes H^1(C_j)\otimes H^1(C_k)\right]_{\lambda_1}
= V_{ij}^k \oplus V_{jk}^i \oplus V_{ki}^j.
$$
It is now easy to see that $\gamma$ is injective. Indeed, by
(\ref{formula}), we see that the images of the maps
$$
H^1(C) \stackrel{p_i^\ast}{\to} H^1(\Delta_{ij})
\stackrel{\gamma_{ij}}{\to} H^3(C^n)
$$
are independent copies of $V$, and also, when $k \not\in \{i,j\}$,
that the image of
$$
H^1(C) \stackrel{p_k^\ast}{\to} H^1(\Delta_{ij})
\stackrel{\gamma_{ij}}{\to} H^3(C^n)
$$
is congruent to $V_{ij}^k$ modulo the sum of the $V_b^a$.
\end{proof}

Similarly, one can show that the $r$ Chern classes of the central
extensions
$$
0 \to \Z^r \to \pi_{g,r}^n \to \pi_g^{r+n} \to 1
$$
are linearly independent in $H^2(\pi_g^{r+n})$ as they correspond to
independent copies of the trivial representation in $H^2(F_g^{r+n})$.
It follows that $H^2(F_{g,r}^n)$ is also pure of weight 2.

Assembling all this, we obtain:

\begin{proposition}
For each choice of a base point $[C]$ of $\M_g$
and for all $g\ge 1$ and $n,r \ge 0$, the natural MHS on
$H^1(F_{g,r}^n)$ is pure of weight 1 and that on $H^2(F_{g,r}^n)$
is pure of weight 2. In addition, the cup product
$$
\Lambda^2 H^1(F_{g,r}^n,\Q) \to H^2(F_{g,r}^n)
$$
is surjective. \qed
\end{proposition}

{}From Morgan's Theorem we deduce that $\p_{g,r}^n$ has a quadratic
presentation for all non-negative $g$, $r$ and $n$.\footnote{In the
genus zero  case, it is well known that $H^1$ has weight 2 and $H^2$
weight 4 as the corresponding classifying spaces are complements
of hyperplanes in affine space.}

Our final task is to determine the relations explicitly. First some
notation. The Lie algebra $\p_{g,r}^n$ is a quotient of the free
Lie algebra generated by
$$
H_1(\p_{g,r}^n) \cong H_1(C^{n+r})
\cong \bigoplus_{i=1}^{n+r} H_1(C_i).
$$
We shall think of elements of $H_1(C^{n+r})$ as indeterminates, and write
them as upper case letters. If $X\in H_1(C)$, we shall denote the
corresponding element of $H_1(C_i)$ by $X\sup{i}$. Fix a symplectic
basis $A_1,\dots,A_g,B_1,\dots, B_g$ of $H_1(C)$. Denote the intersection
number of $X$ and $Y \in H_1(C)$ by $(X\cdot Y)$.

\begin{theorem}
For all $g\ge 1$ and all $r,n\ge 0$,
$$
\Gr^W_\dot \p_{g,r}^n \cong \L(H_1(C)^{\oplus(n+r)})/R
$$
where $R$ is the ideal generated by the relations
\begin{xalignat*}{2}
[X\sup{i},Y\sup{j}] & = [X\sup{j},Y\sup{i}] & \text{all $i$ and $j$;} \cr
[X\sup{i},Y\sup{j}] & =
\frac{(X\cdot Y)}{g} \sum_{k=1}^g\,
[\A{i},\B{j}] & \text{all $i$ and $j$;} \cr
\sum_{k=1}^g\, [\A{i},\B{i}] & = %
\frac{1}{g} \sum_{j\neq i} \sum_{k=1}^g\,
[\A{i},\B{j}] & 1 \le i \le n.
\end{xalignat*}
where $X$ and $Y$ are arbitrary elements of $H_1(C)$.
\end{theorem}

Note that the last relation holds only for those factors corresponding
to a marked point, and not those corresponding to a marked tangent vector.

\begin{proof}
If $\g$ is a graded Lie algebra generated in weight
$-1$ and $H_2(\g)$ of weight 2, then we have an exact sequence
$$
0 \to H_2(\g) \stackrel{\text{cup}^\ast}{\longrightarrow}
\Lambda^2 H_1(\g) \stackrel{\text{bracket}}{\longrightarrow}
Gr^W_{-2}\g \to 0,
$$
where the first map is the dual of the cup product.\footnote{There are
many ways to see this --- the easiest being from the standard complex of
Lie algebra cochains. However, the statement holds in greater generality
--- cf.\ \cite{sullivan:les}.} In our case, the natural injection
$$
H^2(\p_{g,r}^n) \to H^2(F_{g,r}^n,\Q)
$$
is an isomorphism because the cup product is surjective. The
coproduct is the obvious inclusion of $H_2(F_{g,r}^n,\Q)$ into
$\Lambda^2 H_1(C^n,\Q)$, and the sequence is a sequence of $Sp_g$
modules:
$$
0 \to H_2(F_{g,r}^n,\Q) \to \Lambda^2 H_1(C^n,\Q) \to \Gr^W_{-2}
\p_{g,r}^n \to 0.
$$
We will consider one weight at a time. Note that the three weights
occurring in $\Lambda^2 H_1(C^n)$ are 0, $\lambda_2$ and $2\lambda_1$
--- the last being the symmetric square of $H_1(C)$ and second being
the primitive part of $H_2(\Jac C)$. We also have the exact sequence
$$
0 \to H_2(F_{g,r}^n) \to H_2(C^n)
\stackrel{\gamma^\ast}{\to} \bigoplus_{i<j}\Q \to 0
$$
of $Sp_g$ modules. The last map is the dual of the Gysin map.
It follows that
$$
H_2(F_{g,r}^n,\Q)_{\lambda} = H_2(C^n,\Q)_{\lambda}
$$
when $\lambda$ is $2\lambda_1$ or $\lambda_2$.

The $2\lambda_1$ isotypical component is spanned by elements of the
form
$$
X\sup{i}\times Y\sup{j} + Y\sup{i}\times X\sup{j}.
$$
This gives the first relation:
\begin{equation}\label{comm}
[X\sup{i},Y\sup{j}] = [X\sup{j},Y\sup{i}].
\end{equation}

Since $V(\lambda_2)$ is the kernel of the symplectic form
$\Lambda^2 V(\lambda_1) \to \Q$,
the $\lambda_2$ isotypical component of $H_2(C^n)$ is spanned by
vectors of the form
$$
X\sup{i}\times Y\sup{j} - Y\sup{i}\times X\sup{j} -
\frac{(X\cdot Y)}{g} \sum_{k=1}^g
\left( \A{i}\times \B{j} - \B{i}\times \A{j} \right) .
$$
This gives relations of the form
$$
[X\sup{i}, Y\sup{j}] + [X\sup{j}, Y\sup{i}] =
\frac{(X\cdot Y)}{g} \sum_{k=1}^g
\left([\A{i}, \B{j}] + [\A{j}, \B{i}] \right)
$$
which simplifies to the second relation after applying (\ref{comm}).

For the time being, we assume that $r=0$.
The trivial isotypical component lies in an exact sequence
$$
0 \to H_2(F_{g,r}^n)^{Sp_g} \to H_2(C^n)^{Sp_g}
\stackrel{\gamma^\ast}{\to} H^{2g-2}(\Delta) \to 0.
$$
The map $\gamma^\ast$ takes $W\in H_2(C^n)$ to the functional
$$
\{\Delta_{ij} \mapsto W\cdot\Delta_{ij}\}.
$$
Note that
$$
H_2(C^n)^{Sp_g} = \bigoplus_{i=1}^n H_2(C_i) \oplus
\bigoplus_{i<j} \left[H_1(C_i)\otimes H_1(C_j)\right]^{Sp_g}.
$$
The first terms has basis the $Z\sup{i}$, where $Z$ denotes the
integral generator of $H_2(C)$. The second term has basis consisting
of the
$$
Q_{ij} :=
\sum_{k=1}^n \left(\A{i}\times \B{j} - \B{i}\times \A{j}\right).
$$
We next determine a basis of $\ker \gamma^\ast$.

Choose $n$ distinct points $u_1,\dots,u_n$ of $C$. For $i<j$, let
$C_{ij}$ be the image of the map $C\hookrightarrow C^n$ defined by
$$
x \mapsto
(u_1,\dots, u_{i-1},x,u_{i+1},\dots,u_{j-1},x,u_{j+1},\dots,u_n).
$$
Denote its homology class by $Z_{ij}$. It is easily seen that
$$
Z_{ij} = Z_i + Z_j + Q_{ij}.
$$
Observe that
$$
Z_{ij}\cdot \Delta_{kl} =
\begin{cases}
0 & \text{$i$, $j$, $k$, $l$ distinct};\cr
1 & \#\{i,j,k,l\} = 3;\cr
2 - 2g & ij = kl.
\end{cases}
$$
The first two assertions are clear, the second follows from the
projection formula applied to the projection of $C^n$ onto $C^2$
along the $i$th and $j$th factors and the fact that the self
intersection of the diagonal in $C^2$ is the Euler number of $C$.
It follows immediately that a basis of
$$
H_2(F_g^n)^{Sp_g} = \ker\{H_2(C^n)^{Sp_g} \to H^2(\Delta)\}
$$
consists of the
$$
Z_i - \frac{1}{2g}\sum_{j\neq i} Q_{ij}.
$$
These give the relations
$$
\sum_{k=1}^g\, [\A{i},\B{i}] =
\frac{1}{2g} \sum_{j\neq i}
\sum_{k=1}^g\, \left([\A{i},\B{j}] + [\A{j},\B{i}]\right)
$$
which becomes the third relation after an application of (\ref{comm}).

Finally, the third relation is dual to the first Chern class of the
pullback of the tangent bundle of $C$ along $p_i : C^n \to C$. It
follows that in the general case, we do not get any relations coming
from the trivial representation associated to an index corresponding
to a tangent vector.
\end{proof}

We conclude this section with a computation of the generating function
of the lower central series of $\pi_{g,r}^n$. This corrects the formula
in \cite{kohno-oda}. (The galois analogue of this corrected formula is
also stated in \cite[(2.14)]{nak-tak-ueno}.)

\begin{theorem}
For all $g\ge 1$ and all $n\ge 1$, we have
$$
\prod_{k=1}^n \left(1 - 2g\,t -(k-2)t^2\right)
= \prod_{l=1}^\infty\left(1 - t^l\right)^{r_l}
$$
where $r_l$ is the rank of the $l$th graded quotient of the lower
central series of $\pi_g^n$.
\end{theorem}

\begin{proof}
Since $F_g^n$ is a smooth variety and a rational $K(\pi,1)$ by
(\ref{kohno-oda}), we can apply (\ref{wt=lcs}) and the formula
\cite[(9.7)]{hain:cycles} to deduce that
$$
W_{F_g^n}(t) = \prod_{l=1}^\infty\left(1 - t^l\right)^{r_l},
$$
where, for a graded (variation of) MHS $H$
$$
W_H(t) = \sum_{k\ge 0} \chi(\Gr^W_k H)\, t^k
$$
and, for an algebraic variety $X$, $W_X(t) = W_{H^\dot(X)}(t)$.
(Here $\chi$ denotes Euler characteristic.) The result now
follows by induction on $n$ using the fact that
$$
W_{C - \{x_1,\dots,x_n\}}(t) = 1 - 2g\, t + (n-1)t^2
$$
and the following lemma which is proved by induction on the length
of the weight filtration of $V$.
\end{proof}

\begin{lemma}
If $\V$ is an admissible unipotent variation of MHS over a smooth
variety $X$, then
$$
W_{H^\dot(X,\V)}(t) = W_X(t)W_\V(t). \qed
$$
\end{lemma}

\section{The General Case}
\label{decorated}

In this section we assemble results from Sections \ref{special} and
\ref{braids2} to obtain a presentation of $\t_{g,r}^n$ for all $g\ge 6$
and all $r$ and $n \ge 0$. Fix a base point
$$
[C;x_1,\dots,x_n;v_1,\dots,v_n]
$$
of $\M_{g,r}^n$ so that $\t_{g,r}^n$ and $\p_{g,r}^n$, etc.\ all have
mixed Hodge structures. The sequence of Lie algebras
$$
0 \to \p_{g,r}^n \to \t_{g,r}^n \to \t_g \to 0
$$
is exact in the category of MHSs, and therefore remains exact after
applying $Gr^W_\dot$:
$$
0 \to \Gr^W_\dot \p_{g,r}^n \to
\Gr^W_\dot\t_{g,r}^n \to \Gr^W_\dot\t_g \to 0.
$$
By (\ref{canon_split}), $\t_{g,r}^n \otimes \C$ is isomorphic to
the completion of $\Gr^W_\dot\t_{g,r}^n\otimes \C$. So, to find a
presentation of $\t_{g,r}^n$, it suffices to find a presentation of
its associated graded. As in \S\ref{special}, there is a lift of the
canonical homomorphism $\L(H_1(\t_g)) \to \Gr^W_\dot \t_g$ to
a homomorphism $\L(H_1(\t_g)) \to \Gr^W_\dot \t_{g,r}^n$. Given
presentations of $\Gr^W_\dot \t_g$ and $\Gr^W_\dot \p_{g,r}^n$,
a presentation of $\Gr^W_\dot \t_{g,r}^n$ is determined by maps
\begin{gather*}
a : H_1(\t_g)\otimes H_1(\p_{g,r}^n) \to \p_{g,r}^n(2)\\
c : \Lambda^2 H_1(\t_g) \to \p_{g,r}^n(2)_{\lambda_2}\\
c_0 : \Lambda^2 H_1(\t_g) \to \t_{g,r}^n(2)^{Sp_g}
\end{gather*}
induced by the bracket.

Observe that the homomorphism
$$
\p_{g,r}^n \to \p_g^{\oplus(n+r)}
$$
induced by the inclusion $F_{g,r}^n(C) \hookrightarrow C^{n+r}$
induces isomorphisms
$$
H_1(\p_{g,r}^n) \cong H_1(\p_g)^{\oplus(n+r)}\text{ and }
\p_{g,r}^n(2)_{\lambda_2} \to \p_g(2)^{\oplus(n+r)}.
$$
By a naturality argument, the map $a$ is easily seen to be the adjoint
of the map
$$
H_1(\t_g) \to \Hom\left(\bigoplus_{j=1}^{r+n} H_1(\p_g),
\bigoplus_{j=1}^{r+n} \p_g(2)_{\lambda_2}\right)
$$
which is the direct sum of $n+r$ copies of the Johnson homomorphism.

The map $c$ is simply the sum over all the marked points and tangent
vectors
$$
\Lambda^2 H_2(\t_g) \to
\bigoplus_{j=1}^{r+n}\p_g(2) \cong \p_{g,r}^n(2)_{\lambda_2}.
$$
of the maps (\ref{bracket}) which is determined in (\ref{bra_const}).

In remains to determine $c_0$. As in the case of $\t_{g,1}$
considered in Section~\ref{special}, there is a canonical
splitting of the sequence
$$
0 \to \p_{g,r}^n(2)^{Sp_g} \to \t_{g,r}^n(2)^{Sp_g}
\to \t_g^{Sp_g}(2) \to 0
$$
from which we obtain a canonical decomposition
$$
\t_{g,r}^n(2)^{Sp_g} = \p_{g,r}^n(2)^{Sp_g} \oplus \Ga.
$$
As in \S\ref{special}, we identify $H_1(\t_g)$ with the subspace of
$\Lambda^3 V$ which is the kernel of the projection (\ref{projn}),
denote by $\bil \blank \blank$ the unique $Sp_g$ invariant bilinear
form on $H_1(\t_g)$ such that
$$
\bil {a_1\wedge a_2 \wedge a_3} {b_1\wedge b_2 \wedge b_3} = 1,
$$
and choose a generator $\gamma$ of $\Ga$ such that if $u,v\in H_1(\t_g)$,
then the invariant part of $[u,v]$ in $\t_g(2)$ is $\bil u v \,\gamma$.

Observe that there is an exact sequence
$$
0 \to \bigoplus_{1\le i<j\le r+n} \Q(1) \oplus \bigoplus_{i=1}^r \Q(1)
\to \p_{g,r}^n(2) \to \bigoplus_{j=1}^{r+n} \p_g(2) \to 0
$$
of $Sp_g$ modules. It follows that
$$
\p_{g,r}^n(2)^{Sp_g} \cong
\bigoplus_{1\le i<j\le r+n} \Q(1) \oplus \bigoplus_{i=1}^r \Q(1).
$$
The terms indexed by $1\le i \le r$ correspond to the $r$ marked
tangent vectors; those indexed by $1\le i < j \le r+n$
to the diagonals $\Delta_{ij}$. It is easy to see that the composition
$$
\Lambda^2 H_1(\t_g) \to \bigoplus_{i=1}^r \Q(1)
$$
of $c_0$ with the projection
$$
\t_{g,r}^n(2)^{Sp_g} \to \bigoplus_{i=1}^r \Q(1)
$$
is the sum of the maps $c_0$ associated to $\t_{g,1}$ computed
in (\ref{triv_cpt}). So it remains to determine the composition
$$
\Lambda^2 H_1(\t_g) \to \bigoplus_{1\le i<j\le r+n} \Q(1)
$$
of $c_0$ with the projection
$$
e^{r+n} : \t_{g,r}^n(2)^{Sp_g} \to \bigoplus_{1\le i<j\le r+n} \Q(1).
$$
To do this, it suffices to compute $e^2$ in the case of $\t_g^2$,
for then the map $e^{r+n}$ is simply the sum of the $e^2$s over all
diagonals.

In order to compute
$$
e^2 : \t_{g,r}^n(2)^{Sp_g} \to \Q(1)
$$
we use the fact that a punctured tubular neighbourhood of the diagonal
$\Delta$ in $C\times C$ is homeomorphic to the frame bundle of the
tangent bundle of $C$. In this way, we obtain homomorphisms
$$
\t_{g,1} \to \t_g^2 \text{ and } \p_{g,1} \to \p_g^2.
$$
In particular, we have a map
\begin{equation}\label{rest}
\p_{g,1}(2)^{Sp_g} \to \p_g^2(2)^{Sp_g}.
\end{equation}
Using the fact that $\Gr^W_\dot \p_{g,1} \to \Gr^W_\dot \p_g^2$
is a homomorphism and that on $H_1$, it is the diagonal map
$V\to V\oplus V$, we see that the map (\ref{rest}) takes the
generator $\sum_k [A_k,B_k]$ of $\p_{g,1}(2)^{Sp_g}$ to
\begin{multline*}
\sum_{k=1}^g\, \left(
[\A 1, \B 1] + [\A 1, \B 2] + [\A 2, \B 1] + [\A 2, \B 2]\right) \\
= 2 \sum_{k=1}^g\,\left( [\A 1, \B 2] + [\A 2, \B 1]\right).
\end{multline*}
It follows that in $\t_g^2$, the map $c_0$ is given by
$$
c_0[u,v] = \bil u v \,\gamma - \frac{12\bil u v}{g(2g+1)}
\sum_{k=1}^g\, \left([\A 1, \B 2] + [\A 2, \B 1]\right).
$$
This completes the determination of $c_0$ in general and, with it,
the descriptions of the $\t_{g,r}^n$.

\section{Applications}
\label{applications}

\subsection{Cup products and Massey products}
\label{cup}

We have shown that for all $g \ge 6$ and all $r,n\ge 0$, the
Lie algebra $\Gr^W_\dot \t_{g,r}^n$ has a presentation with
only quadratic relations. This implies, using the short exact
sequence in  \cite{sullivan:les} for example, that the cup product
$$
\Lambda^2 \Gr^W_\dot \Hcts^1(\t_{g,r}^n)
\to \Gr^W_\dot \Hcts^2(\t_{g,r}^n)
$$
is surjective. It follows that the cup product
$$
\Lambda^2 \Hcts^1(\t_{g,r}^n) \to \Hcts^2(\t_{g,r}^n)
$$
is also surjective.

Recall that the $l$-fold Massey products constructed from $H^1(A^\dot)$,
where $A^\dot$ is a d.g.a., are defined on a subspace $D_l$ of
$H^1(A^\dot)^{\otimes l}$ and take values in $H^2(A^\dot)/I_{l-1}$,
where $I_{l-1}$ denotes the lift to $H^2(A^\dot)$ of the image of the
Massey products of order $< l$.
It follows that all Massey products in $\Hcts^2(\t_{g,r}^n)$ of order
$\ge 3$ vanish when $g\ge 6$ as the cup product (Massey products of
order 2) map is surjective.

Since the natural map $\Hcts^2(\t_{g,r}^n) \to H^2(T_{g,r}^n,\Q)$
is injective (cf.\ (\ref{cts_ord})) and preserves Massey products,
we have:

\begin{theorem}
For all $g \ge 6$, all Massey products of order $\ge 3$ in
$H^2(T_{g,r},\Q)$ vanish. \qed.
\end{theorem}

\begin{remark}
It follows from the fact that there are non-trivial cubic relations
and no quadratic relations in a minimal presentation of $\t_3$ that
the cup product
$$
\Lambda^2 H^1(T_3,\Q) \to H^2(T_3,\Q)
$$
vanishes, and that the Massey triple product map
$$
H^1(T_3,\Q)^{\otimes 3} \to H^2(T_3,\Q)
$$
is non-trivial.
\end{remark}

It follows from (\ref{lower}) and (\ref{computation}) that for all
$g \ge 6$, $\Hcts^2(\t_g)$ has highest weight decomposition
$$
\Hcts^2(\t_g) \cong V(\lambda_6) + V(\lambda_4) + V(\lambda_2)
+ V(\lambda_2 + \lambda_4).
$$

\begin{theorem}
For all $g \ge 3$, there is an (unnatural) isomorphism
$$
\Hcts^2(\t_{g,r}^n) \cong \Hcts^2(\t_g) \oplus
\left(\Hcts^1(\p_{g,r}^n)\otimes \Hcts^1(\t_g)\right)
\oplus \Hcts^2(\p_{g,r}^n)
$$
of $Sp_g$ modules.
\end{theorem}

\begin{proof}
Chose a base point of $\M_{g,r}^n$. Then
$$
0 \to \p_{g,r}^n \to \t_{g,r}^n \to \t_g \to 0
$$
is an exact sequence of mixed Hodge structures, and the
corresponding spectral sequence
$$
E_2^{s,t} = \Hcts^s(\t_g,\Hcts^t(\p_{g,r}^n)) \implies
\Hcts^{s+t}(\t_{g,r}^n)
$$
is a spectral sequence in the category of mixed Hodge structures.
Since $\t_g$ has negative
weights, the weights on $\Hcts^k(\t_g)$ are $\ge k$. This
and the fact that $H^k(\p_{g,r}^n)$ is a trivial $\t_g$ module
when $k = 1$ and 2 imply that $E_\infty^{s,t} = E_2^{s,t}$ when
$s+t=2$. The result follows.
\end{proof}

\subsection{Period space is not contractible when $g\ge 4$}

Denote by $\h_g$ the Siegel upper half space; that is, the space
of symmetric $g\times g$ complex matrices with positive definite
imaginary part. Denote the image of the period map
$$
\text{Teichm\"uller space } \to \h_g
$$
by $\J_g$, and its closure by $\overline{\J}_g$. When $g \le 3$,
$\Jbar_g=\h_g$, so that $\Jbar_g$ is contractible in these cases.

\begin{theorem}
For each $g\ge 4$, $H^2(\Jbar_g,\Q)$ is non-trivial. Consequently,
$\Jbar_g$ is not contractible.
\end{theorem}

The proof proceeds in two steps. We begin by making a definition.

\begin{definition}
The {\it extended Torelli group} $\Text_g$ is the subgroup of $\Gamma_g$
consisting of those mapping classes which act as $\pm$ the identity
on the first homology of the reference surface.
\end{definition}

We have group extensions
\begin{equation}\label{extensions}
1 \to T_g \to \Text_g \to \Z/2\Z \to 0\text{ and }
1 \to \Text_g \to \Gamma_g \to PSp_g(\Z) \to 1
\end{equation}
where $PSp_g(\Z)$ denotes the quotient of $Sp_g(\Z)$ by $\pm I$.
Note that the first sequence gives rise to a natural action of
$\Z/2\Z$ on $H^\dot(T_g)$.

The first step is:

\begin{proposition}
For all $g \ge 3$, there are natural isomorphisms
$$
H^\dot(\J_g,\Q) \cong H^\dot(\Text_g,\Q) \cong H^\dot(T_g,\Q)^{\Z/2\Z}.
$$
Moreover, when $g\ge 4$, $H^2(\Text_g,\Q)$ is non-trivial.
\end{proposition}

\begin{proof}
Since $g\ge 3$, $\J_g$ is the quotient of Teichm\"uller space by
$\Text_g$. Since the mapping class group acts on Teichm\"uller space
virtually freely, this implies (via standard arguments) that
there is a natural isomorphism
$$
H^\dot(\Text,\Q) \cong H^\dot(\J_g,\Q).
$$
Applying the Hochschild-Serre spectral sequence to the first
of the extensions (\ref{extensions}) above, we see that
$$
H^k(\Text_g,\Q) \cong H^k(T_g,\Q)^{\Z/2}.
$$
But $-I \in Sp_g(\Z)$ acts trivially on $\Hcts^2(\t_g)$, which
implies that
$$
\Hcts^2(\t_g) \subseteq H^2(\Text_g,\Q).
$$
The result follows.
\end{proof}

\begin{remark}
This argument also shows that the image of the cup product
$$
\Lambda^2 H^1(T_3,\Z) \to H^2(T_3,\Z)
$$
is 2-torsion.
\end{remark}

To complete the proof of the theorem, note that $\J_g = \Jbar_g - \cR$
where $\cR$ is the locus of reducible jacobians. By standard arguments,
each component of $\cR$ has complex codimension $\ge 2$ in $\Jbar_g$.
Combining Lefschetz duality and the Gysin sequence, we obtain an exact
sequence
$$
H^{BM}_{6g-k-6}(\cR) \to H^k(\Jbar_g) \to H^k(\J_g)
\to H^{BM}_{6g-k-7}(\cR),
$$
where $H_\dot^{BM}$ denotes Borel-Moore homology. Since $\cR$ has real
codimension 4, it follows that $H^2(\Jbar_g) \cong H^2(\J_g)$. The
theorem follows as $H^2(\J_g)$ is non-trivial.

\subsection{Johnson's conjecture}

In \cite{johnson:survey}, Johnson constructed maps
$$
\phi_k : H_k(T_g) \to H_{k+2}(\Jac S)/[S]\times H_k(\Jac S),
$$
which generalize the classical Johnson homomorphism, which is the
case $k=1$. He conjectured that these homomorphisms are isomorphisms
for all $k$ and sufficiently large $g$.

The following result is an improvement of some unpublished computations
of Morita (cf.\ \cite[\S4]{morita:jap_acad}).

\begin{theorem}
For all $g \ge 3$, the map $\phi_2$ is not injective.
\end{theorem}

\begin{proof}
It is not difficult to see that each $\phi_k$ is $Sp_g(\Z)$
equivariant. Consider its adjoint
$$
\phi_k^t : H^{k+2}(\Jac S,\Q)/\omega\wedge H^k(\Jac S,\Q)
\to H^k(T_g,\Q).
$$
This is also $Sp_g(\Z)$ equivariant. The domain of $\phi_2^t$ is
the primitive cohomology group $PH^4(\Jac S,\Q)$. This is the restriction
to $Sp_g(\Z)$ of the rational representation of $Sp_g$ with highest
weight $\lambda_4$. Since this is an irreducible $Sp_g(\Z)$ module,
the image of $\phi_2^t$ is either isomorphic to $V(\lambda_4)$ or
trivial. But $H^2(T_g,\Q)$ contains the rational representation
$\Hcts^2(\t_g)$. It follows from the results in \S\ref{cup} that
$$
\Hcts^2(\t_g)/\im \phi_2^t \cap \Hcts^2(\t_g)
$$
is non-trivial as it contains
$V(\lambda_6) + V(\lambda_2 + \lambda_4)$ when $g\ge 6$;
$V(\lambda_2 + \lambda_4)$ when $g = 4,5$; and $V(\lambda_3)$
when $g = 3$. The result follows.
\end{proof}

\subsection{Filtrations of $T_g^1$}\label{filtn}

Define a filtration
\begin{equation}\label{filtration}
T_g =
L^1 T_g^1 \supseteq L^2 T_g^1 \supseteq L^3 T_g^1 \supseteq \cdots
\end{equation}
of $T_g^1$ by
$$
L^k T_g^1 =
\left\{\phi \in T_g^1 : \phi_\ast : \pi_1(S,x) \to \pi_1(S,x)
\text{ is congruent to the identity mod } \Gamma^{k+1}\right\}.
$$
It is quite common in the literature for this filtration to be called
the {\it relative weight filtration}, as it is in \cite{asada-nakamura}
and \cite{oda}. In view of (\ref{mhs_torelli}) and (\ref{unequal}), I
feel that this terminology is likely to result in confusion.

\begin{proposition}
This filtration is a descending central series of $T_g^1$ with
torsion free quotients and has the property that
$$
\bigcap_{k=1}^\infty L^k T_g^1
$$
is trivial.
\end{proposition}

\begin{proof}
This follows directly from the fact that the fundamental group of
a compact Riemann surface is residually nilpotent \cite{baumslag},
and the fact that the graded quotients of the lower central series
of a surface group are torsion free \cite{labute}.
\end{proof}

The most rapidly descending series with torsion free quotients
of a group $G$ is the series
$$
G = D^1 G \supseteq D^2 G \supseteq D^3 G \supseteq \cdots
$$
where
$$
D^k G = \{g\in G :
\text{ there is an integer $n>0$ such that } g^n \in \Gamma^k G\}.
$$
This filtration has the property that $D^kG/D^{k+1}$ is the
$k$th term of the lower central series of $G$ mod torsion.
Proofs of these assertions can be found in \cite{passman}.

In the current situation, we have
$$
D^k T_g^1 \subseteq L^k T_g^1.
$$
Johnson's Theorem \cite{johnson:h1} implies that $D^2 T_g^1 =
L^2 T_g^1$ when $g \ge 3$. The computations (\ref{upper}) and
(\ref{lower}) imply that the kernel of $D^2 T_g^1 / D^3 \to
L^2 T_g^1 /L^3$ is isomorphic to $\Z$. Morita was aware of the
fact that the kernel was at least this big --- cf.\ his work on
the Casson invariant \cite{morita:casson}, and asked whether there
is a $k$ such that $D^3 T_g^1 \supseteq L^k T_g^1$. That is, whether
the kernel of $D^2 T_g^1 / D^3 \to L^2 T_g^1 /L^3$ can be detected by
the action of $T_g^1$ on the quotients of $\pi_g$ by the terms of its
lower central series.

More generally, one can ask if the topologies on $T_g^1$ determined by
the filtrations $D^\dot$ and $L^\dot$ are equivalent. (Both are
separated.) That is, for each $k\in \N$, can one find a positive integer
$n(k)$ such that $L^{n(k)} T_g^1 \subseteq D^k T_g^1$~?

Since the groups $T_g^1/D^k$ and $T_g^1/L^k$ are torsion free nilpotent,
they imbed as a Zariski dense subgroup of a unipotent group defined over
$\Q$. One obtains two inverse systems of unipotent groups. It is clear
that the first prounipotent group is the Malcev completion $\T_g^1$ of
$T_g^1$, and the second is the prounipotent group associated to the
pronilpotent Lie algebra $\h_g:=\im\{\t_g^1 \to \d_g\}$, where $\d_g$
is the pronilpotent Lie algebra defined in \S \ref{inf_action}. The two
topologies on $T_g^1$ are equivalent if and only if the natural map
$\t_g^1 \to \h_g$ is an isomorphism. Equivalently, they are equivalent
if and only if $\t_g^1 \to \d_g$ is injective. It is also clear that
the filtration $L^\dot$ of $\t_g^1$ induced from that of $T_g^1$ is
the pullback of the weight filtration of $\h_g^1$, so that
$$
\left(L^kT_g^1/L^{k+1}\right) \otimes\Q \cong \Gr^k_L \t_g^1
\cong \Gr^W_{-k} \h_g^1
$$
and that $L^k\t_g^1 \supseteq \Ga$ for all $k\ge 1$ --- cf.\
(\ref{central_ext}).

\begin{theorem}\label{unequal}
For all $g \ge 3$, and all $k\ge 1$, $L^k\t_g^1 \supseteq \Ga$
so that the natural representation
$\t_g^1 \to \d_g$ is not injective as its kernel contains $\Ga$.
In particular, there is no $k\ge 1$ such that
$W_{-3}\t_g^1 \supseteq L^k\t_g^1$. \qed
\end{theorem}

One can define a filtration $L^\dot$ of $T_g$ by defining $L^k T_g$ to
be the image of $L^k T_g^1$. Using similar arguments, one can prove that
the filtrations $L^\dot$ and $D^\dot$ of $T_g$ do not define equivalent
topologies.

\subsection{A question of Asada and Nakamura}
There is an issue raised by Asada and Nakamura in
\cite[(4.5)]{asada-nakamura} which is closely related to Morita's
question. Denote by $\pi_{g,1}$ the fundamental group $\pi_1(S,v)$ of
$S$ with respect to the tangent vector $v$. It is naturally isomorphic
to the fundamental group of the punctured surface $S$ minus the anchor
point $x$ of $v$. Note that $T_{g,1}$ acts on $\pi_{g,1}$. They define
a filtration $M^\dot$ of $T_g^1$ as follows: First define a filtration
$L^\dot$ of $T_{g,1}$ as in the previous section: $\phi$ is in
$L^kT_{g,1}$ if and only if $\phi$ induces the identity on $\pi_{g,1}$
modulo the $(k+1)$st term of its lower central series. Define $M^kT_g^1$
to be the image of $L^kT_{g,1}$ in $T_g^1$. They then ask whether, after
tensoring with $\Q$, the sequence
$$
0 \to \Gr^W_\dot \pi_g \to \Gr^M_\dot T_g^1
\to \Gr^L_\dot T_g \to 0
$$
is exact. (Recall from (\ref{wt=lcs}) that the lower central series of
$\p_g$ agrees with its weight filtration.) We now give a proof that this
is indeed the case. We continue with the notation of the previous section.

The filtration $M^\dot$ induces a filtration of $\t_g^1$.
Their question then beocmes: is the sequence
$$
0 \to \Gr^W_\dot \p_g \to \Gr^M_\dot \t_g^1
\to \Gr^L_\dot \t_g \to 0
$$
exact? Fix a base point of $\M_{g,1}$ so that $\t_{g,1}$, $\t_g^1$,
$\pi_{g,1}$, $\d_g$, $\p_g$, etc.\ all have compatible MHSs; the MHS
on $\p_{g,1}$ is the limit MHS on $\pi_1(S-\{x\},x_o)$ associated to
the ``degeneration'' where $x_o$ approaches $x$ from the direction of
$v$. Denote the image
of $\t_g^1$ in $\d_g$ by $\h_g^1$, and the image of $\t_g$ in $\o_g$ by
$\h_g$. These have  canonical mixed Hodge structures determined by the
choice of the base point. Since the diagram
$$
\begin{CD}
0 @>>> \p_g @>>> \t_g^1 @>>> \t_g @>>> 0\\
@.	 	@|			@VVV		@VVV \\
0 @>>> \p_g @>>> \h_g^1 @>>> \h_g @>>> 0\\
\end{CD}
$$
commutes and since the top row is exact, it follows that
the bottom row is exact.  Since $\Gr^W_\dot$ is an exact functor,
and since $\Gr^W_k \h_g^n \cong \Gr^L_k\t_g^n$ when
$n = 0,1$, this implies that the sequence
$$
0 \to \Gr^W_\dot \p_g \to \Gr^L_\dot \t_g^1
\to \Gr^L_\dot \t_g \to 0
$$
is exact. To complete the proof, we show that the filtrations $L^\dot$
and  $M^\dot$ of $\t_g^1$ are equal.

Denote by $b_o$ the element of $\pi_{g,1}$ that corresponds to rotating
the tangent vector once about $x$ --- this is a ``Dehn twist about the
boundary of $S-\{x\}$.'' The action of $T_{g,1}$ on $\pi_{g,1}$ fixes
$b_o$, and therefore induces a homomorphism
$T_{g,1} \to \Aut(\pi_{g,1},b_o)$ into the automorphisms of $\pi_{g,1}$
that fix $b_o$. Set $w_o = \log b_o$. This we interpret
as an element of $\p_{g,1}$. The homomorphism above induces
a homomorphism $T_{g,1} \to \Aut(\p_{g,1},w_o)$, and therefore a Lie
algebra homomorphism
$$
\t_{g,1} \to \Der(\p_{g,1},w_o)
$$
into the derivations of $\p_{g,1}$ that annhilate $w_o$. It follows
from standard properties of limit MHSs that $w_o$ spans a copy of
$\Q(1)$ in $\Der\p_{g,1}$. But $\Der(\p_{g,1},w_o)$ is the kernel
of the map $\Der \p_{g,1} \to \p_{g,1}$
that takes $\phi$ to $\phi(w_o) - w_o$. Since $w_o$ is a Hodge class,
this is a morphism of MHS. It follows that $\Der(\p_{g,1},w_o)$ has
a natural MHS.

Since $\p_g$ is the quotient of $\p_{g,1}$ by the ideal generated by
$w_o$ as MHS, the homomorphism
$$
\Der(\p_{g,1},w_o) \to \Der \p_g
$$
is a morphism of MHS. The filtration $M^\dot$ of $\t_{g,1}$ is the
inverse image of the weight filtration under the homomorphism
$\t_{g,1}\to \Der(\p_{g,1},w_o)$. The equality of the filtrations
$L^\dot$ and $M^\dot$ of $\t_g^1$ now follows from the strictness
properties of the weight filtration as the diagram
$$
\begin{CD}
\t_{g,1} @>>> \t_g^1 \cr
@VVV @VVV \cr
\Der(\p_{g,1},w_o) @>>> \Der \p_g \cr
\end{CD}
$$
commutes and all arrows are morphisms of MHS.

\subsection{Cohomology of $\t_g$ and vanishing differentials}

Since $\t_{g,r}^n = W_{-1}\t_{g,r}^n$, it follows that
$$
W_{k-1}\Hcts^k(\t_{g,r}^n) = 0
$$
for all $k\ge 0$. The {\it lowest weight subring of}
$\Hcts^\dot(\t_{g,r}^n)$ is defined to be the subring
$$
\bigoplus_{k\ge 0} W_k\Hcts^k(\t_{g,r}^n).
$$
By \cite[(9.2)]{hain:cycles}, this is a quadratic algebra generated
by $\Hcts^1(\t_{g,r}^n)$ and where the relations are dual to the
second graded quotient of the lower central series of $T_{g,r}^n$.

The following result is a refinement of (\ref{vanishing}).

\begin{theorem}\label{van_diffls}
For each irreducible representation $V(\lambda)$ of $Sp_g$, the
image of the natural homomorphism
$$
\left[W_k\Hcts^k(\t_{g,r}^n) \otimes V(\lambda)\right]^{Sp_g} \to
H^0(Sp_g(\Z),H^k(T_{g,r}^n)\otimes V(\lambda)) = E_2^{0,k}
$$
is contained in
$$
E_\infty^{0,k} = \im\left\{
H^k(\Gamma_{g,r}^n,V(\lambda)) \to H^k(T_{g,r}^n)\otimes V(\lambda)
\right\}.
$$
\end{theorem}

\begin{proof}
Fix a base point of $\M_{g,r}^n$ so that $\t_{g,r}^n$, $\u_{g,r}^n$, etc.\
all have compatible MHSs. Since the extension
$$
0 \to \Ga \to \t_{g,r}^n \to \u_{g,r}^n \to 0
$$
is central with kernel isomorphic to $\Q(1)$, it follows from the
Gysin sequence that the induced map
$$
\bigoplus_{k\ge 0} W_k\Hcts^k(\u_{g,r}^n) \to
\bigoplus_{k\ge 0} W_k\Hcts^k(\t_{g,r}^n).
$$
is surjective, with kernel the ideal generated by the cohomology class
in $W_2\Hcts^2(\u_{g,r}^n)$ corresponding to the extension above.

By (\ref{morph_u}), there is a canonical map
$$
\left[\Hcts^k(\u_{g,r}^n)\otimes V(\lambda)\right]^{Sp_g}
\longrightarrow H^k(\Gamma_{g,r}^n,\V(\lambda)).
$$
The result follows because the diagram
$$
\begin{CD}
\left[W_k\Hcts^k(\u_{g,r}^n)\otimes V(\lambda)\right]^{Sp_g}
@>>> H^\dot(\Gamma_{g,r}^n,\V(\lambda)) \\
@VVV @VVV \\
\left[W_k\Hcts^k(\t_{g,r}^n)\otimes V(\lambda)\right]^{Sp_g}
@>>> H^0(Sp_g(\Z),H^k(T_{g,r}^n)\otimes V(\lambda))
\end{CD}
$$
commutes, and because the left hand vertical map is surjective.
\end{proof}

\subsection{Morita's Conjecture}
We now prove a result which is, in some sense, an affirmation of Morita's
conjecture \cite[2.7]{morita:conj}. Our result is an analogue of his
theorem \cite[6.2]{morita:conj} which is a solution to the conjecture in
the first non-trivial case. He also informs me that he has proved the
second non-trivial case of the conjecture over ${\frac{1}{24}}\Z$.

Suppose that $g\ge 3$. Denote the $k$th term of the lower central
series of $\pi_g$ by $\pi^{(k)}$. Set
$$
\pi_{(k)} = \pi_g/\pi^{(k+1)}.
$$
We know from Labute's theorem \cite{labute} that this is a torsion free
nilpotent group. For each $k\ge 1$, there is a natural representation
$$
\rho_k : \Gamma_g^1 \to \Aut \pi_{(k)}.
$$
The first is simply the standard representation $\Gamma_g^1 \to Sp_g(\Z)$.
Denote the $\Q$ form of the unipotent completion of $\pi_{(k)}$ by $\cP_{(k)}$.
Since $\pi_{(k)}$ is torsion free, the canonical map $\pi_{(k)} \to \cP_{(k)}$
is injective. By the universal mapping property of unipotent completion,
we see that each $\rho_k$ extends to a homomorphism
$$
\rhotilde_k : \Gamma_g^1 \to \Aut \cP_{(k)}.
$$
Denote the Lie algebra of $\cP_{(k)}$ by $\p_{(k)}$. Then
$\Aut \cP_{(k)} \cong \Aut \p_{(k)}$.
It follows that $\Aut \cP_{(k)}$ is a linear algebraic group. Denote
the Zariski closure of the image of $\rhotilde_k$ in this by $G_k$.
It is easy to see that $G_k$ is an extension of $Sp_g(\Q)$ by a
unipotent group:
$$
1 \to U_k \to G_k \to Sp_g(\Q) \to 1
$$
This extension is split exact, so that
$$
G_k \cong Sp_g(\Q) \ltimes U_k.
$$
By the universal mapping property of the relative completion of
$\Gamma_g^1$, there is a homomorphism $\G_g^1 \to G_k$ which
commutes with the projections to $Sp_g$. The following result follows
directly from the fact (\ref{central_ext}) that the natural map
$\T_g^1 \to \U_g^1$ is surjective.\footnote{A direct proof of the
lemma can be given --- cf.\ the proof of \cite[(4.6)]{hain:comp}.}

\begin{lemma}
For each $k\ge 2$, $\rhotilde(T_g^1)$ is Zariski dense in $U_k$. \qed
\end{lemma}

\begin{proposition}
For each $k \ge 2$, the image of $\rhotilde_k$ is a discrete subgroup of
$G_k(\R)$, and the quotient $\im\rho_k\backslash G_k(\R)$ has finite
volume with respect to any left invariant metric on $G_k(\R)$.
\end{proposition}

\begin{proof}
Since every finitely generated subgroup of the $\Q$ points of a unipotent
group $U$ is discrete in $U(\R)$, it follows that $\rhotilde_k(T_g^1)$
is a discrete subgroup of $U_k(\R)$. Since it is also Zariski dense, it
is cocompact. The result now follows as the image of $\Gamma_g^1$ in
$Sp_g(\R)$ is $Sp_g(\Z)$, which is discrete and of finite covolume.
\end{proof}

We should note that Morita works with $\Gamma_{g,1}$ rather than with
$\Gamma_g^1$ as we do. Our arguments work equally well in his case;
we chose to work with $\Gamma_g^1$ as it seems more natural.
In conclusion, we remark that the Lie algebra of $U_k$ is simply the
image $\h_g^1/W_{-k-1}$ of $\u_g^1$ in $\Der \p_{(k)}$. It follows
that the Lie algebra of $U_k$ has a MHS, and is therefore isomorphic
to its associated graded after tensoring with $\C$.

\section{The Universal Connection}
\label{connection}

In this section we construct a universal connection form
$$
\omegatilde \in
E^\dot(\text{Torelli space})\comptensor \Gr^\dot \t_{g,r}^n
$$
with ``scalar curvature''
on Torelli space when $g\ge 3$. Here $E^\dot(X)$ denotes the $C^\infty$
de~Rham complex of a smooth manifold $X$, and $\comptensor$ the completed
tensor product.%
\footnote{The completed tensor product $E^\dot(X)\comptensor\Gr^\dot\g$
is defined to be the inverse limit
$$
\lim_\leftarrow E^\dot(X)\otimes \Gr^\dot\g /\oplus_{l\ge m} \Gr^l\g.
$$
}
This is the analogue of the universal connection
$$
\sum_{ij} d\log(z_i - z_j)\, X_{ij} \in E^\dot(X_n) \otimes \Gr^\dot\p_n
$$
for the braid group $P_n$. Here $X_n$ denotes the classifying space
$$
\C^n - \left\{\text{fat diagonal}\right\}
$$
of the pure braid group, $(z_1,\dots,z_n)$ its coordinates, and $\p_n$ the
Malcev Lie algebra associated to $P_n$. A reasonably precise dictionary
between the case of braid groups and the absolute mapping class groups
$\Gamma_g$ is given in the table.

\begin{table}
{\footnotesize
\begin{tabular}{c|c|p{1.2in}}
\hline
Braid Groups & Mapping Class Groups & \quad Comments \\
 \hline
&&\\
$B_n$ & $\Gamma_g$ & {\rr the group of interest}\\
&&\\
$\Sigma_n$ & $Sp_g(\C)$ & {\rr a semi-simple algebraic group $G$}\\
&&\\
$\rho:B_n \to \Sigma_n$ & $\rho:\Gamma_g \to Sp_g(\C)$ & {\rr homomorphism
to $G$ with dense image}\\
&&\\
$\Sigma_n$ & $Sp_g(\Z)$ & {\rr the image of $\rho$, an arithmetic group}\\
&&\\
$P_n$ & $T_n$ & {\rr the kernel of $\rho$, a residually torsion free
nilpotent group}\\
&&\\
$\cP_n$ & $\U_g$ & {\rr prounipotent radical of the relative completion}\\
&&\\
$B_n \to \Sigma_n \ltimes \cP_n$ &
$\Gamma_g \to \G_g \cong Sp_g(\C)\ltimes \U_g$ &
{\rr the relative completion}\\
&&\\
$\cP_n$ & $\T_g$ & {\rr the unipotent completion of the kernel of $\rho$}\\
&&\\
$\id: \cP_n \to \cP_n$ & $\T_g \to \U_g$ & {\rr the homomorphism to the
prounipotent radical}\\
&&\\
$\p_n$ & $\t_g$ & {\rr the pronilpotent Lie algebra corresponding to
$\ker\rho$}\\
&&\\
$\Gr^\dot \p_n = \L(H_1(P_n))/R$ &
$\Gr^\dot \t_g = \L(H_1(T_g))/R$ &
{\rr quadratic presentations as graded Lie algebras in the category of
representations of $G$}\\
&&\\
$\left\{[X_{ij},X_{kl}],[X_{ij},X_{ik}+X_{jk}]\right\}$ &
$V(\lambda_6), V(\lambda_4), V(\lambda_2), V(\lambda_2 + \lambda_4)$ &
{\rr the quadratic relations}\\
&&\\
$X_n := \C^n - \{\text{fat diagonal}\}$ & $\H_g := \text{ Torelli space}$ &
{\rr the classifying space of the kernel of $\ker\rho$}\\
&&\\
$Y_n := \Sigma_n\backslash X_n$ & $\M_g = Sp_g(\Z)\backslash \H_g$ &
{\rr the classifying space of the group of interest}\\
&&\\
$\sum_{ij} w_{ij}\, X_{ij}
\in E^\dot(X_n)\otimes \Gr^\dot\p_n$ &
$\omega \in E^\dot(\H_g)\comptensor \Gr^\dot\t_g$ &
{\rr the ``universal (projectively) flat connection'' on the classifying
space of $\ker\rho$}
\end{tabular}
}
\end{table}

\begin{question}
The Lie algebra $\Gr^\dot \p_n$ has interesting finite dimensional
representations; namely those associated to Hecke algebras.
Are there analogous representations of $\Gr^\dot \t_{g,r}^n$ where the
canonical central $\Ga$ acts via scalar transformations? These should
lead to interesting projective representations of $\Gamma_{g,r}^n$.
\end{question}

We now give the construction of the connection. First recall that
Torelli space $\H_{g,r}^n$ is the quotient of the Teichm\"uller
space associated to $\Gamma_{g,r}^n$ by $T_{g,r}^n$. It is the
moduli space of isomorphism classes of $(n+r+2g+1)$-tuples
$$
(C;x_1,\dots,x_n;v_1,\dots,v_r;a_1,\dots,a_g,b_1,\dots,b_g)
$$
where $C$ is a compact Riemann surface of genus $g$, $x_1,\dots,x_n$
are $n$ marked points, $v_1,\dots,v_r$ are $r$ marked tangent vectors,
and $a_1,\dots, b_g$ is a symplectic basis of $H_1(C,\Z)$. Since
$T_{g,r}^n$ is torsion free and Teichm\"uller space is contractible,
$\H_{g,r}^n$ is the classifying space of $T_{g,r}^n$.

The bulk of the work needed for the construction of the connection
has already been done in \cite[\S14]{hain:derham}.
Fix a point $x \in \M_{g,r}^n$. It follows from
\cite[\S14.2]{hain:derham} that there is a 1-form
$$
\omega \in E^\dot(\H_{g,r}^n)\comptensor \Gr^W_\dot \u_{g,r}^n
$$
which is integrable and is $Sp_g(\Z)$ invariant in that
\begin{equation}\label{invariance}
s^\ast \omega = Ad(s)\, \omega
\end{equation}
for all $s\in Sp_g(\Z)$. That is, if
$$
\omega = \sum_I w_I X_I, \text{ where } w_I \in E^1(\H_{g,r}^n)
\text{ and } X_I \in \Gr^W_\dot \u_{g,r}^n,
$$
then for all $s\in Sp_g(\Z)$,
$$
\sum_I (s^\ast w_I) X_I = \sum_I w_I (X_I\cdot s^{-1})
$$
where $X\cdot s$ denotes the canonical action of $s\in S$ on
$X\in \u_{g,r}^n$.
This should be compared with the case of braids where the corresponding
formula is easily verified --- cf. \cite[(14.6)]{hain:derham}.

Since there is a canonical isomorphism
$$
\u_{g,r}^n(x) \cong \prod_{l\ge 1}\Gr^W_{-l}\u_{g,r}^n(x)
$$
this form gives rise to a flat connection on the trivial right
$\U_{g,r}^n$ principal bundle
$$
\H_{g,r}^n \times \U_{g,r}^n \to \H_{g,r}^n.
$$
Note that $Sp_g(\Z)$ acts on this bundle via the diagonal action.
The composite
$$
\U_{g,r}^n \to \H_{g,r}^n \to \M_{g,r}^n
$$
is a left principal $Sp_g(\Z)\ltimes \U_{g,r}^n$ bundle (in the
orbifold sense.) The invariance condition (\ref{invariance}) means
that the connection defined by $\omega$ is invariant under the
$Sp_g(\Z)\ltimes \U_{g,r}^n$ action. The monodromy yields a
representation
$$
\Gamma_{g,r}^n \to Sp_g(\C)\ltimes \U_{g,r}^n(x)
$$
As proved in \cite[\S14.2]{hain:derham}, this is the $\C$ form of
the completion of $\Gamma_{g,r}^n$ with respect to the canonical
homomorphism $\Gamma_{g,r}^n \to Sp_g(\C)$.

Since the sequence
$$
0 \to \Ga \to \t_{g,r}^n \to \u_{g,r}^n \to 0
$$
splits canonically over $\C$ (given the choice of the base point
$x$), $\omega$ has a canonical lift
$$
\omegatilde \in E^\dot(\H_{g,r}^n)\comptensor \Gr^W_\dot \t_{g,r}^n
$$
to $\Gr^W_\dot\t_{g,r}^n$.
This form is not integrable, but since $\omega$ is integrable, the
curvature of $\omegatilde$ takes values in the central $\Ga$. It also
has the invariance property (\ref{invariance}).

We will say that a representation
$\phi : \Gr^W_\dot\t_{g,r}^n \to \End(V)$
is {\it projective} if the image of $\Ga$ consists of scalar
matrices. If $V$ is an $Sp_g$ module and $\phi$ is $Sp_g$
equivariant, then $\phi$ should integrate to a homomorphism
$$
\Gamma_{g,r}^n \to PGL(V),
$$
at least when $\phi$ is ``sufficiently small,'' since, in this case,
the composite
$$
\omega_\phi \in E^\dot(X_n)\otimes \End(V)/\text{scalars},
$$
an infinite sum, should converge to an integrable 1-form. The equivariance
of $\phi$ implies that $\omega_\phi$ has the invariance property
(\ref{invariance}), leading to a projectively flat bundle over $\M_{g,r}^n$
with fiber $V$ over the base point $x$.

\appendix

\section{Index of Principal Notation}

This is an index of principal notation. Some notational conventions
appear after the index.

\begin{tabbing}
soenotation \= an explanation ---- quite looooong xxxxxxxxxxxxxxxxxxxx then
\= page reference \kill

$\Gamma_{g,r}^n$ \> mapping class group --- genus $g$, $n$ points, $r$
tangents \> p.~\pageref{group_def}\\

$T_{g,r}^n$ \> Torelli group, genus $g$, $n$ points, $r$ tangents
\> p.~\pageref{torelli_def}\\

$\L(V)$ \> free Lie algebra generated by a vector space $V$ \>
p.~\pageref{lie_def}\\

$F_{g,r}^n$ \> configuration space of points and tangents on a
genus $g$ surface \> p.~\pageref{config_def}\\

$\pi_{g,r}^n$ \> $\pi_1(F_{g,r}^n,\ast)$ --- pure braid group, genus $g$
\> p.~\pageref{fund_def}\\

$\Gamma_{g,r}^n[l]$ \> level $l$ subgroup of $\Gamma_{g,r}^n$ \>
p.~\pageref{level_def}\\

$\G_{g,r}^n$ \> relative completion of $\Gamma_{g,r}^n$ \>
p.~\pageref{comp_def}\\

$\U_{g,r}^n$ \> prounipotent radical of $\G_{g,r}^n$ \>
p.~\pageref{ugp_def}\\

$\u_{g,r}^n$ \> Lie algebra of $\U_{g,r}^n$ \> p.~\pageref{ulie_def}\\

$\T_{g,r}^n$ \> unipotent completion of $T_{g,r}^n$ \>
p.~\pageref{comptor_def}\\

$\t_{g,r}^n$ \> Lie algebra of $\T_{g,r}^n$ \>
p.~\pageref{lietor_def}\\

$\cP_{g,r}^n$ \> unipotent completion of $\pi_{g,r}^n$ \>
p.~\pageref{comppi_def}\\

$\p_{g,r}^n$ \> Lie algebra of $\cP_{g,r}^n$ \>
p.~\pageref{liepi_def}\\

$\M_{g,r}^n[l]$ \> moduli space of decorated genus $g$ curves and
level $l$ structure \> p.~\pageref{mod_def}\\


$\Hcts^\dot(\pi,\Q)$ \> continuous cohomology of a group $\pi$ \>
p.~\pageref{ctsgp_def}\\

$\Hcts^\dot(\g)$ \> continuous cohomology of a Lie algebra $\g$ \>
p.~\pageref{ctslie_def}\\

$\sp_g$ \> symplectic Lie algebra of rank $g$ \> p.~\pageref{symp_def}\\

$R(\sp_g)$ \> representation ring of $\sp_g$ \> p.~\pageref{rep_def}\\

$\lambda_j$ \> $j$th fundamental weight of $\sp_g$ \>
p.~\pageref{wt_def}\\

$V(\lambda)$ \> $\sp_g$ module with highest weight $\lambda$ \>
p.~\pageref{module_def}\\

$|\lambda|$ \> the ``size'' of $V(\lambda)$ \> p.~\pageref{size_def}\\

$V_\lambda$ \> the $\lambda$-isotypical part of an $\sp_g$ module $V$\>
p.~\pageref{iso_def}\\

$\p_g$ \> shorthand for $\p_g^1$ \> p.~\pageref{p_def}\\

$\pi_g$ \> shorthand for $\pi_g^1$ \> p.~\pageref{p_def}\\

$\g(l)$ \> $l$th graded quotient of the weight filtration of $\g$ \>
p.~\pageref{gr_def}\\


$\d_g$ \> the image of $\t_g^1$ in $\Der \p_g$ \> p.~\pageref{der_def}\\

$\o_g$ \> quotient of $\d_g$ by the inner derivations \>
p.~\pageref{out_def}\\

\end{tabbing}

In notation of the form $Y_{g,r}^n$, the decorations $r$ and $n$ are
omitted when they are zero. So, for example, $F_{g,0}^1$ is written
$F_g^1$. In the case of $\pi_{g,r}^n$ and $\p_{g,r}^n$ this is carried
one step further; we denote $\pi_g^1$ by $\pi_g$, and $\p_g^1$ by $\p_g$.
In notation of the form $Y[l]$, the level $l$ is omitted when it is one.
So, for example, $\M_{g,0}^2(1)$ is written $\M_g^2$.

\end{document}